\documentclass[11pt]{article}
\usepackage{amsfonts}

\usepackage{amsfonts,amssymb,amsmath,mathrsfs,multirow,booktabs}

\usepackage{color,latexsym,amsfonts}
 \newcommand{\qed}{\hfill\rule{2mm}{3mm}\vspace{4mm}}

 \newtheorem{theorem}{Theorem}[section]
 \newtheorem{lemma}[theorem]{Lemma}
 \newtheorem{corollary}[theorem]{Corollary}
 \newtheorem{proposition}[theorem]{Proposition}
 \newtheorem{example}[theorem]{Example}
 \newtheorem{Definition}[theorem]{Definition}
 \newtheorem{remark}[theorem]{Remark}
 \newtheorem{condition}[theorem]{Condition}
 \newtheorem{conjecture}[theorem]{Conjecture}

 \def\blemma{\begin{lemma}\sl{}\def\elemma{\end{lemma}}}
 \def\btheorem{\begin{theorem}\sl{}\def\etheorem{\end{theorem}}}
 \def\bcorollary{\begin{corollary}\sl{}\def\ecorollary{\end{corollary}}}
 \def\bdefinition{\begin{Definition}\sl{}\def\edefinition{\end{Definition}}}
 \def\bproposition{\begin{proposition}\sl{}\def\eproposition{\end{proposition}}}
 \def\bremark{\begin{remark}\sl{}\def\eremark{\end{remark}}}
 \def\bcondition{\begin{condition}\sl{}\def\econdition{\end{condition}}}
 \def\bexample{\begin{example}\rm{}\def\eexample{\end{example}}}
 \def\bconjecture{\begin{conjecture}\rm{}\def\econjecture{\end{conjecture}}}

 \def\beqlb{\begin{eqnarray}}\def\eeqlb{\end{eqnarray}}
 \def\beqnn{\begin{eqnarray*}}\def\eeqnn{\end{eqnarray*}}

 \def\mbb{\mathbb}\def\mbf{\mathbf}

 \def\<{\langle}\def\>{\rangle}

 \def\ar{&\!\!}

 \def\eqref#1{{\rm(\ref{#1})}}

 \def\proof{\noindent{\it
 Proof.~}}\def\qed{\hfill$\Box$\medskip}

\def\e{{\mbox{\rm e}}}
\def\<{\left<}\def\>{\right>}
\def\la{\lambda}
 
 \def\mbb{\mathbb}
  \def\mbf{\mathbf}
\newcommand{\dd}{\mathrm{d}}

\newfam\msbmfam\font\tenmsbm=msbm10\textfont
\msbmfam=\tenmsbm\font\sevenmsbm=msbm7
\scriptfont\msbmfam=\sevenmsbm

\def\al{{\alpha}}

\def\la{{\lambda}}

\def\<{\left<}\def\>{\right>}
\def\({\left(}\def\){\right)}

\begin{document}

\centerline{\Large\bf On the extinction-extinguishing dichotomy for a stochastic}

\smallskip\smallskip

\centerline{\Large\bf Lotka-Volterra type population dynamical system}

\bigskip

\centerline{
%Yanxia Ren
Yan-Xia Ren\footnote{LMAM School of Mathematical Sciences \& Center for Statistical Science, Peking University, Beijing, China.
Supported by
NSFC (Nos.~11671017 and 11731009) and LMEQF. Email: yxren@math.pku.edu.cn.},
Jie Xiong\footnote{Department of Mathematics and SUSTech International center for Mathematics, Southern University of Science \& Technology, Shenzhen, China. Southern University of Science and Technology Start up fund Y01286120 and NSFC (Nos.~61873325 and 11831010).  Email: xiongj@sustech.edu.cn},
Xu Yang\footnote{School of Mathematics and Information Science,
North Minzu University, Yinchuan, China. Supported by
NSFC (Nos.~11771018 and 12061004), NSF of Ningxia (No.~2021AAC02018), Major research project for North Minzu University
(No. ZDZX201902) and the Fundamental Research Funds for the Central Universities of North Minzu University (No. 2020KYQD17). Email: xuyang@mail.bnu.edu.cn. Corresponding author.} and Xiaowen Zhou\footnote{Department of Mathematics and Statistics, Concordia
University, Montreal, Canada.
Supported by NSERC (RGPIN-2016-06704).
Email: xiaowen.zhou@concordia.ca.}}

\bigskip\bigskip

{\narrower{\narrower

\noindent{\bf Abstract.}
Applying the Foster-Lyapunov type criteria and a martingale method,
we study a two-dimensional process $(X, Y)$ arising as the unique nonnegative solution to a pair of stochastic differential equations driven by
independent Brownian motions and compensated spectrally positive L\'evy random measures. Both processes $X$ and $Y$ can be identified as
continuous-state nonlinear branching processes where the evolution of $Y$ is negatively affected by $X$. Assuming that process $X$ extinguishes, i.e. it
converges to $0$ but never reaches $0$ in finite time, and process $Y$ converges to $0$,  we identify rather sharp conditions under which the process
$Y$ exhibits, respectively, one of the following behaviors: extinction with probability one, extinguishing with probability one or both extinction
and extinguishing occurring with strictly positive probabilities.

\bigskip


\bigskip

\par}\par}

\section{Introduction and main results}
\setcounter{equation}{0}

\subsection{Introduction on the background, the model and the approach}

Lotka-Volterra model serves as a fundamental ecological system.
The deterministic Lotka-Volterra model for population dynamics describes the evolution of two species suffering from both self-regulations and
interspecific competitions   for limiting resources.
A stochastic Lotka-Volterra process  generalizes the deterministic Lotka-Volterra population dynamics to incorporate the influence of demographic
stochasticity  or  random environmental fluctuations.  In Cattiaux-M\'el\'eard \cite{CaMe10}, an interacting logistic Feller diffusion system is proposed as a
stochastic Lotka-Volterra dynamics whose quasi-stationary distribution is studied. Two different spatial Lotka-Volterra type models are formulated in
Blath et al. \cite{BEM07} as lattice-indexed interacting Feller diffusions and lattice-indexed interacting Fisher-Wright diffusions, respectively, where the
persistence and long term coexistence of the populations are investigated.
 Evans et al. \cite{EHS} consider a two-dimensional diffusion that solves a system of stochastic differential equations with Lotka-Volterra type drift and
 linear diffusion coefficients  driven by a correlated two-dimensional Brownian motion, and study its stationary distribution. Hening and Nguyen \cite{HN18}
 further generalize  the model  of Evans et al. \cite{EHS} and prove results on the rate of convergence to the  stationary distribution.
Similar models have also been studied systematically as solution to a system of stochastic differential equations driven by both Brownian motions and
Poisson random measures. We refer to Zhu and Yin \cite{ZhuYin09} and Bao et al. \cite{BMYY11} and references therein for previous work.

In the above mentioned models, the drift coefficients and (or) the diffusion coefficients are assumed to be of  particular forms.
Hening et al. \cite{HNY18} recently proposed populations dynamics described by $n$-dimensional Kolmogorov systems with nonlinear interactions and driven by
white noise.
Sharp conditions are found for the populations to converge exponentially fast to their stationary distributions and for the populations to converge to
$0$ exponentially fast.
We refer to Bena\"{i}m \cite{Benalm} for a comprehensive study on stochastic persistence and related topics for general interacting SDE systems.

On the other hand, progress has been made on the study of continuous-state branching processes with generalized branching mechanism. The extinction,
explosion and coming down from infinity results for such processes are obtained in Li et al. \cite{LYZh} via martingale approaches. This motivates us to
further study similar behaviors for  the general continuous-state branching processes  with interaction.

In this paper we consider a generalized version of the stochastic competitive Lotka-Volterra process $(X, Y)$ arising as the non-negative, spectrally
positive  solution to a system of stochastic differential equations (SDEs for short) driven by independent Brownian motions and compensated Poisson
random measures.

Intuitively, the process $X$ represents the (re-scaled) size of a population with a certain type of individuals  whose evolution is described by a
continuous-state branching process with a general nonadditive branching mechanism that has been studied in Li et al. \cite{LYZh}.
We also refer to Li \cite{Li_Z2012} for a review on continuous-state branching processes.
Process $Y$ represents a population of another type that is a continuous-state branching process  experiencing a competition pressure  from $X$. From
another point of view, one can also identify $X$ as the environment that affects the evolution of process $Y$.

Process $(X,Y)$ can also be treated as a generalized  two-type
continuous-state nonlinear branching process.
The readers are referred to Ma \cite{Ma2013} and Barczy  et al. \cite{BarczyLi} for two-type
continuous-state branching processes,
and to Li \cite{Li_Z1992}, Hong and Li \cite{HongLi},
Chapter 6 of Li \cite{Li_Z2011} and the
references therein for two-type measure-valued branching processes.

In the study of the Lotka-Volterra process, people are often interested in whether the two different populations still coexist in the long run, or
whether there is only a mono-type population left eventually.
 For a continuous-state branching process people also want to distinguish between extinction and extinguishing that are two distinct ways of converging
 to $0$ as time goes to infinity.
 We say extinction occurs if the process reaches $0$ in finite time, and extinguishing occurs if the process converges to $0$ but never reaches $0$ in
 finite time. In this paper we want to carry out more detailed analysis of the extinction-extinguishing behaviors for process $Y$ given that it
 converges to $0$ eventually, and want to understand how the processes $X$ and $Y$ jointly affect the extinction-extinguishing behaviors of process $Y$.

Note that in  SDE terminology   the above mentioned extinction and extinguishing behaviors correspond  to the accessibility/inaccessibility of boundary $0$ for the  associated SDE. To our best knowledge, such boundary  classifications are rarely known for an interacting system of SDEs with jumps.

 As a first attempt  of studying such interacting population dynamics under general setting, we first consider two populations that both undergo
 nonlinear subcritical branching. We further assume that the interaction between the two populations is one-sided, i.e. the evolution of process $Y$ is
 affected by process $X$ while the impact of $Y$ on $X$ is negligible.
We thus propose and study the following  SDE system:
\begin{equation}\label{1.1}
  \left\{
   \begin{aligned}
   X_t &=X_0 -\int_0^t a_1(X_s) \dd s+\int_0^t a_2(X_s)^{1/2}\dd B_s
   +\int_0^t\int_0^\infty\int_0^{a_3(X_{s-})}z\tilde{M}(\dd s,\dd z,\dd u),\\
   Y_t &=Y_0-\int_0^t[b_1(Y_s)+\theta(Y_s)\kappa(X_s)]\dd s
   +\int_0^t b_2(Y_s)^{1/2}\dd W_s\\ &\quad
   +\int_0^t\int_0^\infty\int_0^{b_3(Y_{s-})}z\tilde{N}(\dd s,\dd z,\dd u),
   \end{aligned}
   \right.
  \end{equation}
where functions  $a_i,b_i$ $(i=1,2,3)$ and $\theta,\kappa$
are nonnegative functions on $[0,\infty)$, $(B_t)_{t\ge0}$ and $(W_t)_{t\ge0}$ are Brownian motions,
$\{\tilde{M}(\dd t,\dd z , \dd u)\}$ and
$\{\tilde{N}(\dd t,\dd z , \dd u)\}$
are compensated Poisson random measures with intensity $\dd t\mu(\dd z)\dd u$
and $\dd t\nu(\dd z)\dd u$, respectively, and with the $\sigma$-finite
non-zero measures $\mu$ and $\nu$ satisfying
 \beqnn
\int_0^\infty(z\wedge z^2)\mu(\dd z)
+\int_0^\infty(z\wedge z^2)\nu(\dd z)<\infty.
 \eeqnn
We also assume that $(B_t)_{t\ge0}$,
$(W_t)_{t\ge0}$,
$\{\tilde{M}(\dd t,\dd z , \dd u)\}$ and
$\{\tilde{N}(\dd t,\dd z , \dd u)\}$ are independent of each other.

Since \eqref{1.1} represents a stochastic continuous-state Lotka-Volterra population system,
by a solution $(X, Y)$ to \eqref{1.1} we mean a  c\`adl\`ag $\mathbb{R}_+^2$-valued process $(X,Y)$ that
satisfies equation \eqref{1.1}  up to the minimum of
the first time of either hitting zero or explosion for both processes $X$ and $Y$, which is a variation of the usual definition of solution to SDE;  see Definition \ref{def}. Conditions on the existence and uniqueness of the solution to
(\ref{1.1}) will be given in Lemma \ref{t5.1}.
Since we are only interested in the solution up to the first time of hitting $0$, the uniqueness holds under mild conditions. The uniqueness of such a solution for SDE had been studied before in Dawson et al. \cite{DawsonLiWang03}
and Li \cite{LiPeisen19}.

The  extinction/extinguishing behaviors of the continuous-state nonlinear branching process $X$  have been studied in Li et al. \cite{LYZh} using a martingale approach.
By imposing conditions on SDEs (\ref{1.1})  so that the solution $X$ extinguishes with probability one and the solution $Y$ converges to $0$ in
probability as time goes to infinity, in this paper we find conditions  under which the process $Y$  becomes extinct in finite time with probability one
and zero, respectively. We further show that under certain conditions, both extinction and extinguishing can happen for $Y$ each with a strictly
positive probability, which is a remarkable phenomena.

For stable Poisson random measures with stable indices in $(1, 2)$ and for power function coefficients in the SDEs in (\ref{1.1}), the conditions can be made
more explicit in terms of the powers and the stable indices, and they turn out to be quite sharp. We are not aware of similar previous results on
solutions to such a  system of general SDEs with jumps.

Our main approach is different from that in Li et al. \cite{LYZh}.
To prove the above mentioned results we first develop stochastic Foster-Lyapunov type  criteria with localized conditions  for  probability of finiteness of the
first time of hitting $0$ by either process $X$ or $Y$. These criteria can be compared with those in
Li et al. \cite{LYZh} for solution to one-dimensional SDE
and are of independent interest.
We refer to Chen \cite{Chen04} and
Meyn and Tweedie \cite{Meyn93} for the (deterministic) Foster-Lyapunov type criteria
for explosion and stability of Markov chains.
The proofs of most of the main results then boil down to finding appropriate  test functions in order to apply the stochastic Foster-Lyapunov type criteria,
 and the localized conditions in the Foster-Lyapunov  criteria  make it more convenient to construct the test functions.

It is remarkable that for the model in Li et al. \cite{LYZh} the Foster-Lyapunov  criteria also produce very sharp results; see recent work in Ma et al. \cite{MYZ20}.

We apply the Foster-Lyapunov criteria to show most of the main results. The key  is to identify the right test functions for which our approach is mostly ad hoc. We  typically start with elementary functions such as power functions and exponential functions, then modify  and (or)  combine these functions in different ways to develop the sharpest possible results. Verification
	  of the criteria often involve lengthy computations.

Among the main results,  applying the stochastic Foster-Lyapunov  criteria we  identify sufficient conditions for the process $Y$  to become extinguishing with probability one or to become extinct with a strictly positive probability.

To find conditions under which the process $Y$ extinguishes with a strictly positive probability, we adopt a different approach, where by first
obtaining an estimate on the time dependent lower bound of the sample pathes of $X$, we apply a martingale argument similar to that in Li et al. \cite{LYZh}
together with a comparison theorem. We also use  either the stochastic Foster-Lyapunov  criteria or the martingale method to study the extinction-extinguishing
behaviors for some critical cases.

The rest of the paper is arranged as follows.  We first present the main results together with
an example of SDEs with power coefficients and stable Poisson random measures in Section 1.2.
The Foster-Lyapunov type criteria are proved in Section 2.
 Proofs of the main results are deferred to Section 3.

\subsection{Main results}

We first present some notations and assumptions.
By Taylor's formula (see \eqref{7.2} and \eqref{7.1} in Section 3 of the following), for $u,z>0$ and $\delta\ge-1$,
 \beqnn
(1+z)^{-\delta}-1+\delta z =
\delta(\delta+1)z^2\int_0^1(1+zv)^{-\delta-2}(1-v)\dd v
 \eeqnn
and
 \beqnn
\ln (1+z)-z=\ln (1+z)-\ln1-z=-z^2\int_0^1(1+zv)^{-2}(1-v)\dd v.
 \eeqnn
We also want to introduce  several auxiliary functions.
For $\delta\in(-1,0)\cup(0,\infty)$, and $u>0$ define
 \begin{gather}
H_{1,\delta}(u)
:=
\frac{1}{\delta(\delta+1)}\int_0^\infty
[(1+z  u^{-1})^{-\delta}-1+\delta z u^{-1}]\mu(\dd z) \nonumber\\
\qquad\qquad=
u^{-2}\int_0^\infty z^2\mu(\dd z)\int_0^1(1+z u^{-1}v)^{-2-\delta}(1-v)\dd v,\label{0.3} \\
H_{2,\delta}(u)
:=
\frac{1}{\delta(\delta+1)}\int_0^\infty [(1+z  u^{-1})^{-\delta}-1+\delta z u^{-1}]\nu(\dd z) \nonumber\\
\qquad\qquad=
u^{-2}\int_0^\infty z^2\nu(\dd z)\int_0^1(1+z u^{-1}v)^{-2-\delta}(1-v)\dd v.\label{0.4}
 \end{gather}
For $u>0$ let
 \begin{gather}
H_{1,0}(u)
 :=-\int_0^\infty \Big(\ln(1+zu^{-1})-zu^{-1}   \Big) \mu(\dd z)\nonumber\\
\qquad\qquad\qquad\quad=
u^{-2}\int_0^\infty z^2\mu(\dd z)\int_0^1(1+z u^{-1}v)^{-2}(1-v)\dd v,\label{0.5}\\
H_{2,0}(u)
 :=-\int_0^\infty \Big(\ln(1+zu^{-1})-zu^{-1}   \Big) \nu(\dd z)\nonumber\\
\qquad\qquad\qquad\quad=
u^{-2}\int_0^\infty z^2\nu(\dd z)\int_0^1(1+z u^{-1}v)^{-2}(1-v)\dd v\label{0.6}
 \end{gather}
and
 \begin{gather}
G_{1,0}(u)
:=
a_1(u)u^{-1}+2^{-1}a_2(u)u^{-2}+a_3(u)H_{1,0}(u),\label{0.7}\\
G_{2,0}(u)
:=
b_1(u)u^{-1}+2^{-1}b_2(u)u^{-2}
+
b_3(u)H_{2,0}(u).\label{0.8}
 \end{gather}
These six functions will appear repeatedly throughout the paper.
The functions $H_{1,\delta}$ and $H_{2,\delta}$
are the same to the function $\delta(\delta+1)H_{\delta-1}$ defined in
(2.1) of \cite{LYZh} which result from Ito's formula applied to power function of $X$.
The functions $G_{1,0}$ and $G_{2,0}$ can be regarded as
the limits of $(1-a)^{-1}G_a$ when $a\to 1$ in (2.3) of \cite{LYZh}. They are also associated to Ito's formula applied to logarithm function of $X$; see \cite{MYZ20}.
To study the extinction-extinguishing phenomena of $Y$ we
impose some conditions on $G_{1,0},G_{2,0}$; see Condition \ref{c} below.

Let $C^2((0,\infty))$
be the space of twice continuously differentiable functions on
$(0,\infty)$
and $C^2((0,\infty)\times(0,\infty))$ denote space of
functions on $(0,\infty)\times(0,\infty)$ with continuous second partial derivatives.

For any generic stochastic  process $Z:=(Z(t))_{t\ge0}$ and constant $w>0$, let
 \beqlb\label{1.2}
\tau_0^Z:=\inf\{t\ge0:Z(t)=0\},\quad
\tau_w^Z:=\tau^Z(w):=\inf\{t\ge0:Z(t)<w\}
 \eeqlb
and
 \beqlb\label{1.3}
\sigma_w^Z:=\sigma^Z(w):=\inf\{t\ge0:Z(t)>w\}
 \eeqlb
with the convention $\inf\emptyset=\infty$.
In the following we state the definition of solution
to SDE \eqref{1.1}, which is defined before the minimum of
the first time of either hitting zero or explosion for the two processes.

\bdefinition\label{def}
By a solution to SDE \eqref{1.1} we mean that a two-dimensional
c\`adl\`ag
process $(X_t,Y_t)_{t\ge0}$ satisfies
SDE \eqref{1.1} up to $\gamma_n:=\tau_{1/n}^X\wedge\tau_{1/n}^Y
\wedge\sigma_n^X\wedge\sigma_n^Y$ for each $n\ge1$ and
$X_t=\limsup_{n\to\infty}X_{\gamma_n-}$
and $Y_t=\limsup_{n\to\infty}Y_{\gamma_n-}$
for $t\ge\lim_{n\to\infty}\gamma_n$.
\edefinition

\begin{remark}\label{r1.1}
	The above definition of solution to SDE \eqref{1.1} allows weaker conditions for uniqueness of solution.  In particular,  the pathwise uniqueness holds if the functions $a_i, b_i, \theta$ and $\kappa$ are all locally Lipschitz on $(0, \infty)$; see Lemma \ref{t5.1}. Also observe that
 \beqnn
\tau_0^X\wedge\tau_0^Y=\lim_{n\to\infty}\tau_{1/n}^X\wedge\tau_{1/n}^Y.
 \eeqnn
\end{remark}

Throughout this paper we assume that the  c\`adl\`ag $\mathbb{R}_+^2$-valued process $(X,Y)$ is the unique solution to  \eqref{1.1},
and consequently, the process $(X,Y)$ has
the strong Markov property.
We always assume that
$X_0,Y_0>0$
and that all the stochastic processes are defined
on the same filtered probability space $(\Omega,\mathscr{F},\mathscr{F}_t,\mbf{P})$.
Let $\mbf{E}$ be the corresponding expectation.

Throughout the paper we also assume that the following conditions hold.
\begin{itemize}
\item[{\normalfont(C1)}]
The functions $a_i,b_i$ $(i=1,2,3)$, $\theta$
and $\kappa$ are nonnegative and bounded on any bounded interval;
\item[{\normalfont(C2)}]
For each $c'>0$,
 \beqnn
\sup_{0<u\le c'}[G_{1,0}(u)+G_{2,0}(u)]<\infty;
 \eeqnn
\item[{\normalfont(C3)}]
For each $0<c'<c''$,
 \beqnn
\inf_{x\in[c',c'']}\{a_2(x)+a_3(x)\}>0,
\inf_{x\in[c',c'']}\{b_2(x)+b_3(x)\}>0~~\mbox{and}~
\mu((c',c''))>0.
 \eeqnn
\end{itemize}

\begin{remark}
Under the above conditions,
with probability one, both processes $X_t$ and $Y_t$
converge to $0$ as $t\to\infty$. But $X_t$ does not become extinct almost surely by \cite[Theorem 2.3 (i) and Proposition 2.6]{LYZh};
see also Lemma \ref{t2.1} in the following.
In this situation we say extinguishing occurs for process $X$.
Note that process $Y$ also becomes extinguishing under the above conditions if $\kappa\equiv 0$. The following theorems give the conditions on the extinction-extinguishing
phenomena of $Y$.
\end{remark}

We first find conditions distinguishing between extinction with probability $0$ and extinction with a positive probability for process $Y$.

\btheorem\label{t0.04}
If there exists a constant
$c^*>0$ so that
$\sup_{0<u\le c^*} \theta(u)u^{-1}<\infty$,
then
$\mbf{P}\{\tau_0^Y<\infty\}=0$.
\etheorem

\btheorem\label{t0.03}
Suppose that there exist constants
$c^*,c_1>0$, $\theta\in[0,1)$ and $\delta>1$ so that
 \beqnn
\inf_{c_1\le u\le c^*}\kappa(u)>0,
~ \inf_{0<u\le c^*}\theta(u)u^{-\theta}>0,\mbox{~and~}
\inf_{0<u\le c^*}\Big[a_2(u)u^{-2-\delta}+a_3(u)u^{-\delta-1}\Big]>0.
 \eeqnn
Then
$\mbf{P}\{\tau_0^Y<\infty\}>0.$
\etheorem

From the above two theorems we find that
the extinction
of $Y$ is caused by $X$ through  the negative drift coefficient function $-\theta(v)\kappa(u)$ for $u$
near zero, and  not caused by the Brownian driven or Poisson-random-measure driven components of the SDE for $Y$. Intuitively, process $Y$ becomes  extinguishing or extinct depending on whether $\theta(u)$ converges to $0$ fast enough or slow enough as $u\to 0+$. Note that under the conditions on $\theta$, the role of function $\kappa$ in these theorems is not essential.

To further study the extinction-extinguishing
behaviors of process $Y$ we need to introduce more sets of conditions.
Since both $X$ and $Y$ have no negative jumps, we  impose upper and lower power function bounds on functions $G_{i,0}(u)$, $\theta(u)$ and $\kappa(u)$ only for $u$ close to $0$. These conditions help to simplify  arguments in proofs and allow more transparent conditions  (in terms of powers of the power functions) for the extinction-extinguishing behaviors for $Y$.
\bcondition\label{c}
\begin{itemize}
\item[{\normalfont(i)}]
There exist constants $\theta\in[0,1)$, $c^*,c_\theta,a,b,\kappa>0$  and $p,q\ge0$ so that
\begin{itemize}
\item[{\normalfont(ia)}]
$G_{1,0}(u)\le a u^p$ for all $0<u\le c^*$;
\item[{\normalfont(ib)}]
$G_{2,0}(u)\ge b u^q$ for all $0<u\le c^*$;
\item[{\normalfont(ic)}]
$\theta(u)\ge c_\theta u^{\theta}$ and $\kappa(u)\ge u^{\kappa}$
for all $0<u\le c^*$.
\end{itemize}
\item[{\normalfont(ii)}]
There exist constants $\theta\in[0,1)$, $c^*,c_\theta,a,b,\kappa>0$  and $p,q\ge0$ so that
\begin{itemize}
\item[{\normalfont(iia)}]
$G_{1,0}(u)\ge a u^p$ for all $0<u\le c^*$;
\item[{\normalfont(iib)}]
$G_{2,0}(u)\le b u^q$ for all $0<u\le c^*$;
\item[{\normalfont(iic)}]
$\theta(u)\le c_\theta u^{\theta}$ and $\kappa(u)\le u^{\kappa}$
for all $0<u\le c^*$.
\end{itemize}
 \item[{\normalfont(iii)}]
Assume that the mapping
$u\mapsto b_3(u)$ is nondecreasing
and that the functions $\theta,b_1,b_2,b_3$ are locally Lipschitz, that is, for each
closed interval $[u,v]\subset (0,\infty)$,
there is a constant $C(u,v)\ge0$ so that
 \beqnn
|\theta(x)-\theta(y)|+\sum_{i=1,2,3}|b_i(x)-b_i(y)|\le C(u,v)|x-y|
 \eeqnn
for all $u\le x,y\le v$.
\end{itemize}
\econdition

We remark that the above Condition \ref{c} (iii) is needed for a  comparison theorem, Proposition \ref{t4.4}, which is applied in proofs for
Theorems \ref{t0.02} and \ref{t0.05}.

The following theorems further distinguish between extinction with a positive probability and extinction with probability one for $Y$. In particular, we identify conditions under which both extinction and extinguishing happen with a strictly positive probability.

\btheorem\label{t0.01}
Suppose
that Condition \ref{c} (i) holds with
 \beqlb\label{1.5}
\frac{q\kappa}{q+1-\theta}<p.
 \eeqlb
Then
$\mbf{P}\{\tau_0^Y<\infty\}=1.$
\etheorem

\btheorem\label{t0.02}
Suppose that Condition \ref{c} (ii) and (iii) hold
with
 \beqlb\label{1.4}
\frac{q\kappa}{q+1-\theta}>p>0.
 \eeqlb
Then
$\mbf{P}\{\tau_0^Y<\infty\}<1.$
\etheorem

In the following we consider Condition \ref{c} for either
$pq=0$ or $\frac{q\kappa}{q+1-\theta}=p$ with $p,q>0$.
Observe that the case for $p>0$
and $q=0$ is addressed in Theorem \ref{t0.01} on the extinction behavior.

\btheorem\label{t0.06}
Suppose that Condition \ref{c} (i) holds for constants satisfying  one of the following conditions:
\begin{itemize}
\item[{\normalfont(i)}]
$p=q=0$ and $b/a>\kappa/(1-\theta)$,
\item[{\normalfont(ii)}] $p,q>0$, $\frac{q\kappa}{q+1-\theta}=p$, and
 \beqlb\label{1.6}
\frac{ap}{q(q+1-\theta)}
<\Big(\frac{b}{1-\theta}\Big)^{\frac{1-\theta}{q+1-\theta}}
\cdot
\Big(\frac{c_\theta }{q}\Big)^{\frac{q}{q+1-\theta}}.
 \eeqlb
\end{itemize}
Then
$\mbf{P}\{\tau_0^Y<\infty\}=1$.
\etheorem

\btheorem\label{t0.05}
Suppose that Condition \ref{c} (ii) and (iii) hold
for constants satisfying  one of the following conditions:
\begin{itemize}
\item[{\normalfont(i)}]
$p=q=0$ and $b/a<\kappa/(1-\theta)$,
\item[{\normalfont(ii)}] $p=0$, $q>0$.
\end{itemize}
Then
$\mbf{P}\{\tau_0^Y<\infty\}<1$.
\etheorem

 Given the above theorems, there are still cases for the parameters $a, b, c_\theta, p, q, \theta, \kappa $ in which the extinction-extinguishing behaviors are unknown.

\bconjecture\label{conje}
We have
$\mbf{P}\{\tau_0^Y<\infty\}<1$  if Condition \ref{c} (ii) and (iii) hold with $p,q>0$, $\frac{q\kappa}{q+1-\theta}=p$ and
 \beqnn
\frac{ap}{q(q+1-\theta)}
>\Big(\frac{b}{1-\theta}\Big)^{\frac{1-\theta}{q+1-\theta}}
\cdot
\Big(\frac{c_\theta}{q}\Big)^{\frac{q}{q+1-\theta}}.
 \eeqnn
\econjecture

To better understand the conditions, we next consider an example of
SDE system (\ref{1.1}) with power function coefficients and stable Poisson random measures.
\bexample
Suppose that there are constants $a_i,b_i,\theta\ge0$,
$\kappa,\eta>0$, $\alpha_1,\alpha_2\in(1,2)$ and $q_i,p_i\ge0$ so that
$\kappa(u)=u^\kappa$, $\theta(u)=\eta u^\theta$,
 \beqnn
a_i(u)=a_iu^{p_i+i},\, b_i(u)=b_iu^{q_i+i} \mbox{ for }i=1,2\mbox{ and }
a_3(u)=a_3u^{p_3+\al_1},\,b_3(u)=b_3u^{q_3+\al_2},
 \eeqnn
and
 \beqnn
\mu(\dd z)=\frac{\alpha_1(\alpha_1-1)}{\Gamma(2-\alpha_1)}
z ^{-1-\alpha_1}\dd z ,~~
\nu(\dd z)=\frac{\alpha_2(\alpha_2-1)}{\Gamma(2-\alpha_2)}
z ^{-1-\alpha_2}\dd z
 \eeqnn
with Gamma function $\Gamma$.
We also assume that $a_2+a_3>0$ and $b_2+b_3>0$.
Then $H_{1,0}(u)=\Gamma(\alpha_1)u^{-\alpha_1}, H_{2,0}(u)=\Gamma(\alpha_2)u^{-\alpha_2}$ and
 \beqnn
G_{1,0}(u)
=
a_1u^{p_1}+2^{-1}a_2u^{p_2}+a_3\Gamma(\alpha_1)u^{p_3}, ~~
G_{2,0}(u)
=
b_1u^{q_1}+2^{-1}b_2u^{q_2}
+
b_3\Gamma(\alpha_2)u^{q_3}.
 \eeqnn
Let
 \beqnn
p
 \ar:=\ar
\min\{p_11_{\{a_1\neq0\}},p_21_{\{a_2\neq0\}},p_31_{\{a_3\neq0\}}\}, \cr
q
 \ar:=\ar
\min\{q_11_{\{b_1\neq0\}},q_21_{\{b_2\neq0\}},q_31_{\{b_3\neq0\}}\}
 \eeqnn
and
 \beqnn
a
 \ar:=\ar
a_11_{\{p_1=p\}}+\frac{a_2}{2}1_{\{p_2=p\}}
+a_3\Gamma(\alpha_1)1_{\{p_3=p\}},\cr
b
 \ar:=\ar
b_11_{\{q_1=q\}}+\frac{b_2}{2}1_{\{q_2=q\}}
+b_3\Gamma(\alpha_2)1_{\{q_3=q\}}.
 \eeqnn
Note that  constants $p$ and $q$ are the minimum powers of the power functions
in the expressions of $G_{1,0}$ and $G_{2,0}$, respectively, and the associated power function or functions dominate the behaviours of the corresponding polynomial $G_{i, 0}(u)$  for $u$  near zero.
The constants $a$ and $b$ represent  coefficients
of  the (possibly combined) dominant power functions, respectively.
Since processes $X$ and $Y$ have no negative jumps, these dominant power functions together determine the extinction/extinguishing behaviours of $Y$.

Combining Theorems \ref{t0.04}-\ref{t0.03} and \ref{t0.01}-\ref{t0.05},  we have
\begin{itemize}
\item[{\normalfont(i)}]
$\mbf{P}\{\tau_0^Y<\infty\}=0$ if $\theta\ge1$;
\item[{\normalfont(ii)}]
$\mbf{P}\{\tau_0^Y<\infty\}>0$ if $0\le\theta<1$;
\item[{\normalfont(iii)}]
$\mbf{P}\{\tau_0^Y<\infty\}=1$ if $0\le\theta<1$ and one of the following holds:
\begin{itemize}
\item[(iiia)] $p=q=0$ and $b/a>\kappa/(1-\theta)$;
\item[(iiib)] $p>0$ and $q=0$;
\item[(iiic)]
$p,q>0$ and $\frac{q\kappa}{q+1-\theta}<p$;
\item[(iiid)]
$p,q>0$, $\frac{q\kappa}{q+1-\theta}=p$ and
 \beqnn
\frac{ap}{q(q+1-\theta)}
<\Big(\frac{b}{1-\theta}\Big)^{\frac{1-\theta}{q+1-\theta}}
\cdot
\Big(\frac{\eta }{q}\Big)^{\frac{q}{q+1-\theta}};
 \eeqnn
\end{itemize}
\item[{\normalfont(iv)}]
$0<\mbf{P}\{\tau_0^Y<\infty\}<1$ if $0\le\theta<1$ and one of the following holds:
\begin{itemize}
\item[(iva)] $p=q=0$ and $b/a<\kappa/(1-\theta)$;
\item[(ivb)] $p=0$ and $q>0$;

\item[(ivc)] $p,q>0$ and $\frac{q\kappa}{q+1-\theta}>p$.
\end{itemize}
\end{itemize}
\eexample

 \begin{remark}
From the above example we have the following insights.
The extinction of $Y$ is caused by relatively large negative interaction
$-\eta Y_t^\theta X_t^\kappa$. If $\theta>0$ is small enough,
$Y_t^\theta$ decreases slowly enough as $Y_t\to 0+$ and there
is enough negative drift to cause extinction.
\begin{itemize}
\item
If $\kappa>0$ is further relatively small, then $X_t^\kappa$ also decreases slowly enough
as $X_t\to0+$ so that there is a large enough negative  drift $-\eta Y_t^\theta X_t^\kappa$ that causes extinction for $Y$ with probability one.
\item
On the other hand, if $\kappa>0$ is not relatively small, then the negative drift $-\eta Y_t^\theta X_t^\kappa$  becomes small enough
when $X_t$ starts to take small values. In this case,   $Y$ can survive with a positive probability.
\end{itemize}
\end{remark}

\section{Two-dimensional stochastic Foster-Lyapunov type criteria}
\setcounter{equation}{0}

The study of boundary behaviors for Markov processes started with the boundary classification of Markov chains and	some
Chinese probabilist had made important contributions on it.	Among them
Mu-Fa	Chen identified  the explosion/non-explosion conditions for continuous time Markov chains and Markov jump processes in 1980s; see the review paper Chen \cite{Chen15} and the book Chen \cite{Chen04}  and references therein where the uniqueness and non-uniqueness problems  essentially correspond to the non-explosion and the explosion, respectively, for Markov chains.
	Khasminskii \cite{Khas} proved similar conditions for diffusion processes. 	These conditions were later referred to as Foster and Lyapunov criteria for more general Markov process; see e.g. Meyn and Tweedie  \cite{Meyn93}.

Using
one-dimensional Foster-Lyapunov type criteria,
 an estimate is found in Section 4 of \cite{LYZh}  on the first passage probabilities for the continuous-state nonlinear branching process $X$.
    %probability of
% one-dimensional stochastic process smaller than   strictly positive constant and exceed a level
A  Foster-Lyapunov type criterion is also identified in
 \cite{MYZ20}  for non-extinction of the continuous-state nonlinear branching process.
These criteria generalize similar results for Markov chains;
see Chen \cite[Theorems 2.25 and 2.27]{Chen04}.
The conditions for Propositions \ref{t0.2} and \ref{t0.3} in the following
are also similar to those of \cite[Theorems 2.25 and 2.27]{Chen04}
which are the criteria of uniqueness of $q$-process.

In this section we establish the two-dimensional criteria,
which will be used to prove Theorems \ref{t0.04}, \ref{t0.03}, \ref{t0.01} and \ref{t0.06}.
In the following,
let  $(x_t,y_t)_{t\ge0}$ with $x_0,y_0>0$ denote a two-dimensional  Markov process where
 $(x_t)_{t\ge0}$ and $(y_t)_{t\ge0}$ are  two nonnegative
processes defined before the minimum of their
first times of hitting $0$ or explosion. Let
$L_t$  be an operator such that
for each $g\in C^2((0,\infty)\times(0,\infty))$ and $m,n\ge1$, the process
$t\mapsto M_{t\wedge \gamma_{m,n}}^g$ is a local martingale,
where
 \beqlb\label{2.4}
M_t^g:=
g(x_t,y_t)-g(x_0,y_0)
-\int_0^tL_sg(x_s,y_s)\dd s
 \eeqlb
and $\gamma_{m,n}:=\tau_n\wedge\sigma_m$
with $\tau_n:=\tau_{1/n}^x\wedge\tau_{1/n}^y$
and $\sigma_m:=\sigma_m^x\wedge\sigma_m^y$. Then a natural candidate for $L_t$ is the generator of the process $(x_t, y_t)_{t\ge0}$.
Define the stopping time $\tau_0:=\tau_0^x\wedge\tau_0^y$.
Since the two processes $(x_t)_{t\ge0}$ and $(y_t)_{t\ge0}$ are defined before the first time of hitting zero or explosion,
 \beqlb\label{2.6}
\tau_0=\lim_{n\to\infty} \tau_n.
 \eeqlb

\bproposition\label{t0.2}
Suppose that there is a non-negative function $g\in C^2((0,\infty)\times(0,\infty))$
and a sequence of positive  constants $(d_m)_{m\geq 1}$ satisfying
\begin{itemize}
\item[{\normalfont(i)}]
$\lim_{x\wedge y\to0+}g(x,y)=\infty$;
\item[{\normalfont(ii)}]
$L_tg(x,y)\leq d_mg(x,y)$  for all $t>0,\,x,y\in(0,m)$ and all large $m\geq 1$.
\end{itemize}
Then $\mbf{P}\{\tau_0<\infty\}=0$.
\eproposition
\proof
Observe that there is a sequence of stopping times $(\gamma_k)_{k\ge1}$ so that
$\gamma_k\to\infty$ almost surely as $k\to\infty$ and $t\mapsto M_{t\wedge \gamma_{m,n,k}}^g$
is a martingale for each $m,n,k\ge1$, where $\gamma_{m,n,k}:=\gamma_{m,n}\wedge\gamma_k$.
By \eqref{2.4} and condition (ii), for each $m,n,k\ge1$ and $t\ge0$,
 \beqlb\label{2.5}
\mbf{E}\big[g(x_{t\wedge \gamma_{m,n,k}},y_{t\wedge \gamma_{m,n,k}})\big]
 \ar=\ar
g(x_0,y_0)+ \int^{t}_0 \mbf{E}\big[L_sg(x_s,y_s)
1_{\{s\leq \gamma_{m,n,k}\}}\big]\dd s \cr
 \ar\leq\ar
g(x_0,y_0)+d_m\int^{t}_0
\mbf{E}\big[g(x_s,y_s)
1_{\{s\leq \gamma_{m,n,k}\}}\big]\dd s \cr
 \ar\leq\ar
g(x_0,y_0)+d_m\int^{t}_0 \mbf{E}
\big[g(x_{s\wedge \gamma_{m,n,k}},y_{s\wedge \gamma_{m,n,k}})\big]\dd s.
 \eeqlb
Using \eqref{2.5} and Gronwall's lemma we obtain that for all $k\ge1$,
 \beqnn
\mbf{E}\big[g(x_{t\wedge \gamma_{m,n,k}},y_{t\wedge \gamma_{m,n,k}})\big]
\le
g(x_0,y_0)\e^{d_mt},\qquad t\ge0.
 \eeqnn
Letting $k\to\infty$ we have
 \beqnn
\mbf{E}\big[g(x_{t\wedge \gamma_{m,n}},y_{t\wedge \gamma_{m,n}})\big]
\le
g(x_0,y_0)\e^{d_mt},\qquad t\ge0,
 \eeqnn
which implies  that for each $m\geq 1$,
 \beqlb\label{0.2}
~\mbf{E}\Big[\lim_{n\to\infty}g(x_{t\wedge \gamma_{m,n}},y_{t\wedge \gamma_{m,n}})\Big]
\le
\liminf_{n\to\infty}
\mbf{E}\big[g(x_{t\wedge \gamma_{m,n}},y_{t\wedge \gamma_{m,n}}))\big]
\le g(x_0,y_0)\e^{d_mt}
 \eeqlb
by Fatou's lemma.
From condition (i) and \eqref{2.6} it follows that $\mbf{P}\{\tau_0>t\wedge\sigma_m\}=1$ for each $m\ge1$ and $t>0$.
Letting $t\to\infty$ we get $\mbf{P}\{\tau_0\ge\sigma_m\}=1$ for each $m\ge1$.
Thus, $\tau_0\ge\lim_{m\to\infty}\sigma_m$ almost surely.
Since these two processes are defined before the first time of hitting zero or explosion, then
$\mbf{P}\{\tau_0=\infty$ or $\lim_{n\to\infty}\sigma_n=\infty\}=1$.
This concludes that
$\tau_0=\infty$ almost surely.
\qed

\bproposition\label{t0.3}
Suppose that $\sup_{t\ge0} (x_t + y_t)<\infty$ almost surely.
We also assume that there exist a nonnegative function $g\in C^2((0,\infty)\times(0,\infty))$ and a sequence of nonnegative processes $(d_m)_{m\geq 1}$ satisfying the following conditions:
\begin{itemize}
\item[{\normalfont(i)}] $0<\sup_{x,y>0} g(x,y)<\infty$;

\item[{\normalfont(ii)}] $\int_0^\infty d_m(t)\dd t=\infty$ almost surely for all $m\ge1$;

\item[{\normalfont(iii)}] $L_tg(x_t,y_t)\geq d_m(t)g(x_t,y_t)$ for all
$0<t<\sigma_m$ and large $m\ge1$.
\end{itemize}
Then
$\mbf{P}\{\tau_0<\infty\}\ge g(x_0,y_0)/ \sup_{x,y>0}g(x,y)$.
\eproposition
\proof
Let $D_m(t):=\int_0^t d_m(s)\dd s$.
Then for all $m\ge1$,
 \beqlb\label{2.3}
D_m(t)\to \infty\mbox{  almost surely as } t\to\infty
 \eeqlb
by condition (ii).
Let $([M^g,M^g]_t)_{t\ge0}$ be the quadratic variation process of $(M_t^g)_{t\ge0}$.
Then the mapping $t\mapsto[M^g,M^g]_t$ is right continuous with left limits.
It follows that
 \beqlb\label{2.7}
\sup_{s\in[0,t]}[M^g,M^g]_s<\infty
 \eeqlb
almost surely for all $t>0$.
For $n,m,k\ge1$ define stopping times $\gamma_n$ and $\gamma_{m,n,k}$ by
 \beqnn
\gamma_k:
=\inf\{t\ge0: [M^g,M^g]_t\ge k\},\quad \gamma_{m,n,k}:=\gamma_{m,n}\wedge\gamma_k.
 \eeqnn
Then
 \beqlb\label{2.8}
[M^g,M^g]_t\le k\quad\mbox{ for all }0\le t<\gamma_k\mbox{~~ and ~~}k\ge1
 \eeqlb
and
 \beqlb\label{2.9}
\lim_{k\to\infty}\gamma_k=\infty
 \eeqlb
almost surely by \eqref{2.7}.
One can see that by the assumptions,
 \beqlb\label{2.2}
\lim_{m\to\infty}\mbf{P}\{\sigma_m<\infty\}
\le\lim_{m\to\infty}\mbf{P}\big\{\sup_{t\ge0}(x_t+y_t)\ge m\big\}=0.
 \eeqlb
Moreover, by \cite[p. 73]{P05},
$t\mapsto M_{t\wedge \gamma_{m,n,k}}^g$ is a martingale,
where $M_t^g$ is defined in \eqref{2.4}.
It follows from integration by parts that
 \beqlb\label{2.1}
 \ar\ar
g(x_{t\wedge\gamma_{m,n,k}},
y_{t\wedge\gamma_{m,n,k}})\e^{-D_m(t)} \cr
 \ar=\ar
g(x_0,y_0)+\int_0^tg(x_{s\wedge\gamma_{m,n,k}},
y_{s\wedge\gamma_{m,n,k}})\dd (\e^{-D_m(s)}) \cr
 \ar\ar
+\int_0^t\e^{-D_m(s)}\dd
(g(x_{s\wedge\gamma_{m,n,k}},y_{s\wedge\gamma_{m,n,k}})) \cr
 \ar=\ar
g(x_0,y_0)
-\int_0^tg(x_{s\wedge\gamma_{m,n,k}},
y_{s\wedge\gamma_{m,n,k}})d_m(s)\e^{-D_m(s)}\dd s \cr
 \ar\ar
+\int_0^t\e^{-D_m(s)}L_sg(x_s,y_s)1_{\{s\leq \gamma_{m,n,k}\}}\dd s
+\int_0^t\e^{-D_m(s)}\dd M^g_{s\wedge\gamma_{m,n,k}}.
 \eeqlb
By the Burkholder-Davis-Gundy inequality and
\eqref{2.8}, there is a constant $C>0$ so that
for all $T>0$,
 \beqnn
 \ar\ar
\mbf{E}\Big[\sup_{0\le t\le T}\Big|\int_0^t\e^{-D_m(s)}\dd M^g_{s\wedge\gamma_{m,n,k}}\Big|^2\Big] \cr
 \ar\ar\quad\le
C\mbf{E}\Big[\Big|\int_0^T\e^{-2D_m(s)}
\dd [M^g,M^g]_{s\wedge\gamma_{m,n,k}}\Big|\Big] \cr
 \ar\ar\quad\le
C\mbf{E}\Big[[M^g,M^g]_{T\wedge\gamma_{m,n,k}}\Big]\le Ck.
 \eeqnn
It then follows from \cite[p. 38]{P05} that
$t\mapsto\int_0^t\e^{-D_m(s)}\dd M_{s\wedge\gamma_{m,n,k}}^g$ is a martingale.
Taking expectations on both sides of \eqref{2.1} we get
 \beqnn
 \ar\ar
\int_0^t\mbf{E}\Big[d_m(s)\e^{-D_m(s)}
g(x_{s\wedge\gamma_{m,n,k}},y_{s\wedge\gamma_{m,n,k}})\Big]\dd s
+
\mbf{E}\Big[g(x_{t\wedge\gamma_{m,n,k}},y_{t\wedge\gamma_{m,n,k}})\e^{-D_m(t)}\Big] \cr
 \ar\ar\qquad=
g(x_0,y_0)
+\int_0^t \mbf{E}\Big[\e^{-D_m(s)}L_sg(x_s,y_s)1_{\{s\leq \gamma_{m,n,k}\}}\Big]\dd s.
 \eeqnn
Letting $t\to\infty$ and using condition (i), \eqref{2.3}
and the dominated convergence theorem we get
 \beqnn
 \ar\ar
\int_0^\infty\mbf{E}\Big[d_m(t)\e^{-D_m(t)}
g(x_{t\wedge\gamma_{m,n,k}},y_{t\wedge\gamma_{m,n,k}})\Big]\dd t \cr
 \ar\ar\quad
=
g(x_0,y_0)
+\int_0^\infty \mbf{E}\Big[\e^{-D_m(t)}L_tg(x_t,y_t)1_{\{t\leq \gamma_{m,n,k}\}}\Big]\dd t.
 \eeqnn
Using condition (iii) we have
 \beqnn
 \ar\ar
\int_0^\infty\mbf{E}\Big[d_m(t)\e^{-D_m(t)}
g(x_{t\wedge\gamma_{m,n,k}},y_{t\wedge\gamma_{m,n,k}})\Big]\dd t \cr
 \ar\ar\quad
\ge
g(x_0,y_0)
+\int_0^\infty \mbf{E}\Big[d_m(t)\e^{-D_m(t)}g(x_t,y_t)1_{\{t\leq \gamma_{m,n,k}\}}\Big]\dd t,
 \eeqnn
which implies
 \beqnn
g(x_0,y_0)
 \ar\le\ar
\mbf{E}\Big[\int_0^\infty d_m(t)\e^{-D_m(t)}g(x_{\gamma_{m,n,k}},y_{\gamma_{m,n,k}})
1_{\{t> \gamma_{m,n,k}\}}\dd t\Big] \cr
 \ar\le\ar
c_0\mbf{E}\Big[\int_{\gamma_{m,n,k}}^\infty  d_m(t)\e^{-D_m(t)}\dd t\Big]
=c_0\mbf{E}\big[\e^{-D_m(\gamma_{m,n,k})}\big]
 \eeqnn
by condition (i) and \eqref{2.3} again, where $c_0:=\sup_{x,y>0} g(x,y)$.
Letting $n,k\to\infty$ and using \eqref{2.6} and \eqref{2.9} we get
 \beqnn
g(x_0,y_0)
 \ar\le\ar
c_0\mbf{E}\big[\e^{-D_m(\tau_0\wedge\sigma_m)}\big]
=
c_0\mbf{E}\big[\e^{-D_m(\tau_0\wedge\sigma_m)}
(1_{\{\sigma_m<\infty\}}+1_{\{\sigma_m=\infty\}})\big] \cr
 \ar\le\ar
c_0\mbf{P}\{\sigma_m<\infty\}
+
c_0\mbf{E}\big[\e^{-D_m(\tau_0)}\big].
 \eeqnn
By \eqref{2.2}, for each $\varepsilon\in(0,1)$, there is a large enough $m\ge1$
so that
 \beqnn
c_0\mbf{P}\{\sigma_m<\infty\}\le \varepsilon g(x_0,y_0),
 \eeqnn
which means that
 \beqnn
(1-\varepsilon)g(x_0,y_0)\le c_0\mbf{E}\big[\e^{-D_m(\tau_0)}\big]
\le
c_0\mbf{E}\big[\e^{-D_m(\tau_0)}1_{\{\tau_0=\infty\}}
+1_{\{\tau_0<\infty\}}\big]
=
c_0\mbf{P}\{\tau_0<\infty\},
 \eeqnn
where \eqref{2.3} is used in the last equation.
Taking $\varepsilon\to0$ one ends the proof. \qed

By an  argument similar to that in the proof of Proposition \ref{t0.3},
we obtain the next result.
\bcorollary\label{cor2.1}
Suppose that $\sup_{t\ge0} (x_t+y_t)<\infty$ almost surely
and $g\in C^2((0,\infty)\times(0,\infty))$ is a nonnegative function
with $0<\sup_{x,y>0}g(x,y)<\infty$.
If there exist a constant $\varepsilon>0$ and a nonnegative function $h$ on $(0,\infty)$ so that
 \beqnn
L_tg(x_t,y_t)\geq h(x_t)g(x_t,y_t),\qquad 0<t<\sigma_\varepsilon
 \eeqnn
and
$\int_0^\infty h(x_t\wedge\varepsilon)\dd t=\infty$ almost surely,
then
 \beqnn
\mbf{P}\{\tau_0\wedge\sigma_\varepsilon<\infty\}
\ge g(x_0,y_0)/\sup_{x,y>0}g(x,y).
 \eeqnn
\ecorollary
\proof
We can prove the assertion with
$d_m(t)$ and $\tau_0$ respectively replaced by
 $h(x_t\wedge\varepsilon)$ and $\tau_0\wedge\sigma_\varepsilon$
in the proof of Proposition \ref{t0.3}.
We leave the details of the proof to the readers.
\qed

Similar to Propositions \ref{t0.2} and \ref{t0.3},
we can also obtain the associated   assertions for the one-dimensional processes.
Suppose that $x:=(x_t)_{t\ge0}$ is a non-negative Markov process and the operator $L_t$
is defined in the following: for each $g\in C^2((0,\infty))$ and $m,n\ge1$,
$t\mapsto M_{t\wedge \tilde{\gamma}_{m,n}}^g$ is a local martingale,
where
 \beqnn
M_t^g:=
g(x_t)-g(x_0)
-\int_0^tL_sg(x_s)\dd s
 \eeqnn
and $\tilde{\gamma}_{m,n}:=\tau_{1/n}^x\wedge\sigma_m^x$.

\bcorollary\label{cor2.2}
Suppose that there are a non-negative function $g\in C^2((0,\infty))$
and constants $d_m\ge0$, $m\ge1$  satisfying
$\lim_{y\to0}g(y)=\infty$ and
$L_tg(y)\leq d_mg(y)$ for all $m\ge1$, $y\in(0,m)$ and $t>0$.
Then $\mbf{P}\{\tau_0^x<\infty\}=0$.
\ecorollary

\bcorollary\label{cor2.3}
Suppose that $\sup_{t\ge0} x_t<\infty$ almost surely,
and that there exist a nonnegative function $g\in C^2((0,\infty))$
and a sequence of nonnegative processes $(d_n)_{n\ge1}$
so that $0<\sup_{y>0}g(y)<\infty$,
$\int_0^\infty d_n(t)\dd t=\infty$ almost surely
and $L_tg(x_t)\geq d_n(t)g(x_t)$ for all
$0<t<\sigma_n$ and $n\ge1$.
Then
$\mbf{P}\{\tau_0^x<\infty\}\ge g(x_0)/\sup_{x>0}g(x)$.
\ecorollary

\section{Proofs of the main results}
\setcounter{equation}{0}

In this section we establish the proofs of
Theorems \ref{t0.04}-\ref{t0.03}
and \ref{t0.01}-\ref{t0.05}.
We first state some notations and assertions which will be used in the proofs.
For $g\in C^2((0,\infty)\times(0,\infty))$ we define
 \begin{equation}\label{3.1}
 \begin{split}
&K_z^1g(x,y):=g(x+z,y)-g(x,y)-zg_x'(x,y),\\
&K_z^2g(x,y):=g(x,y+z)-g(x,y)-zg_y'(x,y)
 \end{split}
 \end{equation}
for $x,y,z>0$
and
 \beqlb\label{3.20}
Lg(x,y):=L_1g(x,y)+L_2g(x,y)
 \eeqlb
with
 \beqlb\label{3.21}
L_1g(x,y)
 :=
-a_1(x)g'_{x}(x,y)+\frac12a_2(x)g''_{xx}(x,y)
+a_3(x)\int_0^\infty K^1_zg(x,y)\mu(\dd z)
 \eeqlb
and
 \beqlb\label{3.22}
L_2g(x,y)
 \ar:=\ar
-[b_1(y)+\kappa(x)\theta(y)]g'_{y}(x,y) \cr
 \ar\ar
+\frac12b_2(y)g''_{yy}(x,y)
+b_3(y)\int_0^\infty K^2_zg(x,y)\nu(\dd z),
 \eeqlb
where $g'_x,g_{xx}''$ and $g'_y,g_{yy}''$ denote the
first and the second partial derivatives of $g$ with respect to $x$ and $y$.
By \eqref{1.1} and It\^o's formula, $L$ is the generator of $(X,Y)$
and independent of time $t$.
By Taylor's formula, for any bounded continuous second derivative function $g$,
 \beqlb\label{7.2}
K_zg(x)=z^2\int_0^1g''(x+zv)(1-v)\dd v,
 \eeqlb
where
  \beqlb\label{7.1}
K_zg(x):=g(x+z)-g(x)-zg'(x),\qquad x,z>0.
 \eeqlb
Since for all $x\in\mbb{R}$,
$\e^x-1\ge x$,
then by \eqref{7.2}, for all $x,y,z,\lambda>0$,
 \beqlb\label{7.7}
 \ar\ar
\e^{\lambda y^r}[\e^{-\lambda (y+z)^r}
-\e^{-\lambda y^r}
+\la r zy^{r-1}\e^{-\lambda y^r}]\cr
 \ar\ar\qquad=
[\e^{\lambda y^r-\lambda (y+z)^r}
-1
+\la rz y^{r-1}] \cr
 \ar\ar\qquad\ge
-\lambda [(y+z)^r-y^r-r zy^{r-1}] \cr
 \ar\ar\qquad=
\lambda r(1-r) z^2 \int_0^1(y+zv)^{r-2}(1-v)\dd v \cr
 \ar\ar\qquad=
\lambda r(1-r) z^2 y^{r-2}\int_0^1(1+zy^{-1}v)^{r-2}(1-v)\dd v.
 \eeqlb
Moreover, for $0<r<1$,
 \beqlb\label{7.7b}
 \ar\ar
\e^{\lambda y^r}[\e^{-\lambda (y+z)^r}
-\e^{-\lambda y^r}
+\la r zy^{r-1}\e^{-\lambda y^r}] \cr
 \ar\ar\quad\ge
\lambda r(1-r) z^2 y^{r-2}\int_0^1(1+zy^{-1}v)^{-2}(1-v)\dd v.
 \eeqlb

\subsection{Preliminary Results}

\blemma\label{t2.3}
For any $u,v\ge0$ and $\bar{p},\bar{q}>1$ with $1/\bar{p}+1/\bar{q}=1$, we have
 \beqnn
u+v\ge \bar{p}^{1/\bar{p}}\bar{q}^{1/\bar{q}} u^{1/\bar{p}}v^{1/\bar{q}}.
 \eeqnn
\elemma
\proof
The above inequality  follows from the Young inequality.
\qed

Recall the definitions of $\tau_0^Z$ and $\tau_w^Z$ in \eqref{1.2}
for the process $Z$ and  constant $w>0$.
\blemma\label{t2.1}
\begin{itemize}
\item[{\normalfont(i)}]
For any $0<w_1<X_0$ and $0<w_2<Y_0$, we have
$\tau_{w_1}^X<\infty$ and $\tau_{w_2}^Y<\infty$ almost surely.
Moreover, $\lim_{t\to\infty}X_t=0$ and $\lim_{t\to\infty}Y_t=0$  almost surely.
\item[{\normalfont(ii)}]
$\mbf{P}\{\tau_0^X=\infty\}=1$ and $\mbf{P}\{\tau_0^Y=\infty\}=1$ if $\kappa(x)=0$ for all $x>0$.
\end{itemize}
\elemma
\proof
Observe that for each $w>1$, $(1+zu^{-1}v)^{-1-w}\le(1+zu^{-1}v)^{-2}$.
Then
 \beqnn
 \ar\ar
u^{-2}\int_0^\infty z^2\mu(\dd z)\int_0^1(1+z u^{-1}v)^{-1-w}(1-v)\dd v
\le
H_{1,0}(u), \cr
 \ar\ar
u^{-2}\int_0^\infty z^2\nu(\dd z)\int_0^1(1+z u^{-1}v)^{-1-w}(1-v)\dd v
\le
H_{2,0}(u),
 \eeqnn
where $H_{1,0}$ and $H_{2,0}$ are defined in \eqref{0.5} and \eqref{0.6}.
It is obvious that condition (C2) is stronger than the assumption
in Theorem 2.3 (i) of \cite{LYZh}.
Thus, under conditions (C2) and (C3),
applying Theorem 2.3 (i) and Proposition 2.6 of \cite{LYZh} we prove the assertions.
\qed

The next result gives an estimate on the maximums of processes $X$ and $Y$.

\blemma\label{t3.2}
Given $X_0,Y_0>0$, for any $\delta\in(0,1/2)$ and $\varepsilon>0$,
there exists a constant
$C>0$ that does not depend on $\varepsilon,X_0,Y_0$
so that
 \beqnn
\mbf{P}\Big\{\sup_{t\ge0}X_t\ge\varepsilon\Big\}
\le C\varepsilon^{-\delta}X_0^{\delta},\quad
\mbf{P}\Big\{\sup_{t\ge0}Y_t\ge\varepsilon\Big\}
\le
C\varepsilon^{-\delta}Y_0^{\delta}.
 \eeqnn
\elemma
\proof
Observe that
 \beqnn
(X_s+z )^{2\delta}-X_s^{2\delta}-2\delta X_s^{2\delta-1} z
=X_s^{2\delta}\big[(1+zX_s^{-1})^{2\delta}-1+(-2\delta) z X_s^{-1}\big]
 \eeqnn
and then
 \beqnn
\int_0^\infty
\big[(X_s+z )^{2\delta}-X_s^{2\delta}-2\delta X_s^{2\delta-1} z\big]
\mu(\dd z)
=-2\delta(1-2\delta)X_s^{2\delta}H_{1,-2\delta}(X_s),
 \eeqnn
where the function $H_{1,-2\delta}$ is defined in \eqref{0.3}.
Then by \eqref{1.1} and It\^o's formula, we have
 \beqlb\label{3.7}
X_t^{2\delta}
 \ar=\ar
X_0^{2\delta}-2\delta\int_0^t a_1(X_s)X_s^{2\delta-1}\dd s
-\delta(1-2\delta)\int_0^t a_2(X_s)X_s^{2\delta-2}\dd s \cr
 \ar\ar
-2\delta(1-2\delta)\int_0^ta_3(X_s)X_s^{2\delta}H_{1,-2\delta}(X_s)\dd s
+2\delta\int_0^t a_2(X_s)^{1/2}X_s^{2\delta-1}\dd B_s \cr
 \ar\ar
+\int_0^t\int_0^\infty\int_0^{a_3(X_{s-})}
[(X_{s-}+z )^{2\delta}-X_{s-}^{2\delta}]\tilde{M}(\dd s,\dd z,\dd u)  \cr
 \ar=:\ar
X_0^{2\delta}-\sum_{i=1}^5 A_i(t,2\delta).
 \eeqlb

Since $0<\delta<1/2$ and $a_1,a_2,a_3$ are nonnegative by condition (C1), then $A_i(t,2\delta)\ge0$ for all $t\ge0$ and $i=1,2,3$.
With $2\delta$ replaced by $\delta$ in \eqref{3.7} it follows that
 \beqlb\label{3.12}
X_t^\delta
\le
X_0^\delta+|A_4(t,\delta)|+|A_5(t,\delta)|.
 \eeqlb
For all $n\ge1$ let $\tilde{\tau}_n:=\tau^X_{1/n}\wedge\sigma^X_n$.
Then
 \beqnn
\mbf{E}[A_4(t\wedge\tilde{\tau}_n,2\delta)]
=
\mbf{E}[A_5(t\wedge\tilde{\tau}_n,2\delta)]
=0.
 \eeqnn
 It follows from \eqref{3.7} that
$\mbf{E}[A_2(t\wedge\tilde{\tau}_n,2\delta)]
\le X_0^{2\delta}$.
We then apply
Fatou's lemma to get
 \beqlb\label{3.10}
\mbf{E}\Big[\sup_{t\ge0}A_2(t,2\delta)\Big]
=\mbf{E}\Big[\lim_{t,n\to\infty}A_2(t\wedge\tilde{\tau}_n,2\delta)\Big]
\le
\liminf_{t,n\to\infty}\mbf{E}[A_2(t\wedge\tilde{\tau}_n,2\delta)]
\le X_0^{2\delta}.
 \eeqlb
Similarly, we can also get
 \beqlb\label{3.11}
\mbf{E}\Big[\sup_{t\ge0}A_3(t,2\delta)\Big]
\le X_0^{2\delta}.
 \eeqlb
By the Burkholder-Davis-Gundy inequality, the H\"older inequality and the estimate \eqref{3.10},
there is a constant $C_1>0$ so that
 \beqlb\label{3.8}
\mbf{E}\Big[\sup_{t\ge0}|A_4(t,\delta)|\Big]
 \ar\le\ar
\delta C_1\mbf{E}\Big[\Big|\int_0^\infty a_2(X_s)X_s^{2\delta-2}\dd s\Big|^{1/2}\Big] \cr
 \ar\le\ar
\delta C_1\Big|\mbf{E}\Big[\int_0^\infty a_2(X_s)X_s^{2\delta-2}\dd s\Big]\Big|^{1/2} \cr
 \ar\le\ar
\delta C_1[\delta(1-2\delta)]^{-1/2}\Big|\mbf{E}\Big[\sup_{t\ge0}A_2(t,2\delta)
\Big]\Big|^{1/2} \cr
 \ar\le\ar
\delta C_1[\delta(1-2\delta)]^{-1/2} X_0^\delta.
 \eeqlb

For fixed $z,x\ge0$ and $u\ge0$ let $h(u)=(1+zx^{-1}u)^{2\delta-2}$.
Then $h$ is decreasing and $(h(u)-h(1/2))(1/2-u)\ge0$
for all $u\ge0$. It follows that
 \beqnn
\int_0^1(h(u)-h(1/2))(1/2-u)\dd u\ge0.
 \eeqnn
Since $\int_0^1h(1/2)(1/2-u)\dd u=0$,
 then $\int_0^1h(u)(1/2-u)\dd u\ge0$.
Moreover,
\[\int_0^1h(u)(1-u)\dd u=
\frac12\int_0^1h(u)\dd u+\int_0^1h(u)(1/2-u)\dd u
\ge \frac12\int_0^1h(u)\dd u.\]
Combining this with Taylor's formula and  H\"older's inequality we further have
 \beqnn
 \ar\ar
[(x+z)^\delta-x^\delta]^2
 =
\delta^2\Big|z\int_0^1(x+zu)^{\delta-1}\dd u\Big|^2 \cr
 \ar=\ar
\delta^2x^{2\delta-2}z^2\Big|\int_0^1(1+zx^{-1}u)^{\delta-1}\dd u\Big|^2 \cr
 \ar\le\ar
\delta^2x^{2\delta-2}z^2\int_0^1(1+zx^{-1}u)^{2\delta-2}\dd u \cr
 \ar\le\ar
2\delta^2x^{2\delta-2}z^2\int_0^1(1+zx^{-1}u)^{2\delta-2}(1-u)\dd u.
 \eeqnn
Then by the Burkholder-Davis-Gundy inequality and the H\"older inequality,
there are constants $C_2>0$ and $C_3=C_3(\delta)>0$ so that
 \beqlb\label{3.13}
\mbf{E}\Big[\sup_{t\ge0}|A_5(t,\delta)|\Big]
 \ar\le\ar
\delta
C_2\mbf{E}\Big[\Big|\int_0^\infty a_3(X_s)\dd s\int_0^\infty
[(X_s+z )^{\delta}-X_s^{\delta}]^2\mu(\dd z)\Big|^{1/2}\Big] \cr
 \ar\le\ar
\delta
C_2\Big|\mbf{E}\Big[\int_0^\infty a_3(X_s)\dd s\int_0^\infty
[(X_s+z )^{\delta}-X_s^{\delta}]^2\mu(\dd z)\Big]\Big|^{1/2} \cr
 \ar\le\ar
2^{1/2}\delta^2 C_2
\Big|\mbf{E}\Big[\int_0^\infty a_3(X_s)X_s^{2\delta}H_{1,-2\delta}(X_s)\dd s\Big]\Big|^{1/2} \cr
 \ar=\ar
C_3\Big|\mbf{E}\Big[\sup_{t\ge0}A_3(t,2\delta)\Big]\Big|^{1/2}
\le
C_3 X_0^\delta,
 \eeqlb
where \eqref{3.11} is used in the last inequality.
Combining \eqref{3.13} with \eqref{3.12} and \eqref{3.8}
we have
 \beqnn
\mbf{E}\Big[\sup_{t\ge0}X_t^\delta\Big]
\le C_4X_0^\delta
 \eeqnn
for some constant $C_4>0$ independent of $X_0$.
Then by the Markov inequality,
 \beqnn
\mbf{P}\Big\{\sup_{t\ge0}X_t>\varepsilon\Big\}
\le
\varepsilon^{-\delta}\mbf{E}\Big[\sup_{t\ge0}X_t^\delta\Big]
\le
C_4\varepsilon^{-\delta}X_0^\delta.
 \eeqnn
By the same argument we can show that
 \beqnn
\mbf{P}\Big\{\sup_{t\ge0}Y_t>\varepsilon\Big\}
\le
C_5\varepsilon^{-\delta}Y_0^\delta
 \eeqnn
for some constant $C_5>0$.
This concludes the proof.
\qed

\subsection{Proof of Theorem \ref{t0.04}}

%{\red Our proof is an application of Proposition \ref{t0.2}.}

\noindent{\it Proof of Theorem \ref{t0.04}.}
We apply Proposition \ref{t0.2} to prove Theorem \ref{t0.04}.
The key is to construct a function $g$ that satisfies the conditions
(i) and (ii)
in Proposition \ref{t0.2}.
For $\rho>0$ we choose the function $g$ as
 \beqnn
g(x,y)
=
x^{-\rho}+y^{-\rho}+1,\qquad x,y>0.
 \eeqnn
Then for all $x,y>0$,
 \beqnn
 \begin{aligned}
 \ar\ar
g'_x(x,y)=-\rho x^{-\rho-1},\quad
g''_{xx}(x,y)=\rho(\rho+1) x^{-\rho-2}, \cr
 \ar\ar
g'_y(x,y)=-\rho y^{-\rho-1},\quad~
g''_{yy}(x,y)=\rho(\rho+1)y^{-\rho-2}.
 \end{aligned}
 \eeqnn
It thus follows that for $x,y>0$,
 \beqlb\label{6.11}
-g'_x(x,y)
\le
\rho x^{-1}g(x,y),\quad
g''_{xx}(x,y)
\le
\rho(\rho+1) x^{-2}g(x,y)
 \eeqlb
and
 \beqlb\label{6.12}
- g'_{y}(x,y)
\le
\rho y^{-1}g(x,y), \quad~\,
g''_{yy}(x,y)\le
\rho(\rho+1)y^{-2}g(x,y).
 \eeqlb
Moreover, by \eqref{3.1},
for $x,y,z>0$,
 \beqnn
K_z^1g(x,y)
 \ar=\ar
z^2\int_0^1g_{xx}''(x+z v,y)(1-v)\dd v \cr
 \ar=\ar
\rho(\rho+1) z^2\int_0^1(x+z v)^{-\rho-2}(1-v)\dd v \cr
 \ar\le\ar
\rho(\rho+1) g(x,y)x^{-2}z^2\int_0^1(1+z x^{-1}v)^{-2}(1-v)\dd v.
 \eeqnn
Thus
 \beqlb\label{6.2}
\int_0^\infty K_z^1g(x,y)\mu(\dd z)
\le
\rho(\rho+1)g(x,y)H_{1,0}(x),\quad x,y>0,
 \eeqlb
where the function $H_{1,0}$ is defined in \eqref{0.5}.
Similarly, we can obtain
 \beqlb\label{6.3}
\int_0^\infty K_z^2g(x,y)\nu(\dd z)
\le
\rho(\rho+1)g(x,y)
H_{2,0}(y),\quad x,y>0,
 \eeqlb
where the function $H_{2,0}$ is defined in \eqref{0.6}.

Recalling \eqref{3.20}-\eqref{3.22} and
combining \eqref{6.11}-\eqref{6.3}, we get
 \beqnn
L_1g(x,y)
\le
\rho(\rho+1)G_{1,0}(x)g(x,y)
 \eeqnn
and
 \beqnn
L_2g(x,y)
\le
\rho(\rho+1)G_{2,0}(y)g(x,y)
+\rho\kappa(x)\theta(y)y^{-1} g(x,y)
 \eeqnn
for $x,y>0$,
where $G_{1,0}$ and $G_{2,0}$ are defined in \eqref{0.7} and
\eqref{0.8}, respectively.
Then under conditions (C1) and (C2) and the assumption of the theorem,
$g(x,y)^{-1}Lg(x,y)$ is bounded for $x,y$ in any bounded interval.
Therefore, for each $n\ge1$,
there is a constant $d_n>0$ so that
$Lg(x,y)\le d_ng(x,y)$ for all $x,y\in(0,n)$,
which implies condition (ii) in Proposition \ref{t0.2}.
By the definition of $g$, condition (i) in Proposition \ref{t0.2}
are obvious. Since $\tau_0^X=\infty$, $\mbf{P}$-a.s. by Lemma \ref{t2.1} (ii), we have
$\tau_0=\tau_0^Y$, $\mbf{P}$-a.s.,
where $\tau_0=\tau_0^X\wedge\tau_0^Y$ by Remark \ref{r1.1}.
It follows from the assertion in Proposition \ref{t0.2} that
$\mbf{P}\{\tau_0^Y=\infty\}=1$.
\qed

\subsection{Proof of Theorem \ref{t0.03}}

\noindent{\it Proof of Theorem \ref{t0.03}.}
We want to apply Proposition \ref{t0.3} where  the key is
to construct a function $g$ satisfying the conditions (i)-(iii) in Proposition \ref{t0.3}.
We assume that the constant $c^*$ in the assumption
satisfies $0<c^*<1$.
By condition (C3), there is a small enough constant $c_2\in(0,c^*)$ so that
$\int_{c_2}^{c^*}\mu(\dd z)>0$.
Let $0<c_1<c_2<c_3<c^*$.
Choose a constant $c_0>0$ so that
 \beqlb\label{6.9}
\inf_{c_1\le u\le c^*}\kappa(u)\ge c_0,
\inf_{0<u\le c^*}\theta(u)u^{-\theta}\ge c_0
 \eeqlb
and
 \beqlb\label{6.10}
\inf_{0<u\le c^*}\big[a_2(u)u^{-2-\delta}+a_3(u)u^{-\delta-1}\big]\ge c_0,
 \eeqlb
where $\delta>1$ is the constant appearing in the assumption.
The proof is given in the following three steps.

{\bf Step 1.}
In this step we construct the function $g$ and summarize  some of its properties.
Let $g_0\in C^2((0,c^*))$
satisfy $g_0(x)=x^{-\delta}$ for $x\in (0,c_2)$
and $g_0(x)=(x-c^*)^{-2}$ for $x\in(c_3,c^*)$.
We choose function $g_0$ so that $g_0,g_0'$ and $g_0''$ are all bounded in $[c_2,c_3]$.
For $\lambda_1,\lambda_2>1,~ \bar{c}:=\pi/(2c^*)$   and $0<r<1-\theta$, define a nonnegative function $g$ by
 \beqnn
g(x,y)
:=
\exp\{-\lambda_1 g_0(x)-\lambda_2 (\tan \bar{c}y)^r\}1_{\{x,y<c^*\}},
\qquad
x,y>0,
 \eeqnn
where we only need the properties of a $\tan$ function such
that it is equivalent to $x$ near zero and is infinite at $\pi/2$.
Then $g\in C^2((0,\infty)\times(0,\infty))$,
and for $0<x,y<c^*$,
 \beqlb\label{6.5}
g'_x(x,y)
 =
-\lambda_1g_0'(x) g(x,y),\quad
g'_y(x,y)
 =
-\lambda_2\bar{c}r(\tan \bar{c}y)^{r-1} (\cos \bar{c}y)^{-2} g(x,y),
 \eeqlb
and
 \begin{gather}
g''_{xx}(x,y)
 =
\lambda_1 \big[\lambda_1|g'_0(x)|^2
- g''_0(x)\big]g(x,y),\qquad \qquad ~\label{4.13} \\
g''_{yy}(x,y)
 =
\lambda_2r\bar{c}^2g(x,y)(\sin \bar{c}y)^{2r-2}(\cos \bar{c}y)^{-2-2r}[\lambda_2 r \nonumber\\ \qquad \qquad \qquad \qquad\quad
+(1-r)(\sin \bar{c}y)^{-r}(\cos \bar{c}y)^r-2(\sin \bar{c}y)^{2-r}(\cos \bar{c}y)^r] \nonumber\\\qquad \qquad \qquad
 \ge
\lambda_2r\bar{c}^2g(x,y)(\sin \bar{c}y)^{2r-2}(\cos \bar{c}y)^{-2-2r}(\lambda_2 r-2) \nonumber\\
\ge
2^{-1}(\lambda_2r\bar{c})^2g(x,y)(\sin \bar{c}y)^{2r-2}~~~~ \label{4.7}
 \end{gather}
as $\lambda_2\ge 4r^{-1}$.
Observe that the constant $\delta>1$ by the assumption of the theorem.
Taking $\lambda_1$ large enough so that
$\lambda_1\delta
- (\delta+1)c_2^{\delta}
\ge\lambda_1$ and $2\lambda_1\ge3|c^*|^2$ in the following,
by \eqref{4.13} we get
 \beqlb\label{4.5}
g''_{xx}(x,y)/g(x,y)
 \ar=\ar
\lambda_1^2\delta^2 x^{-2\delta-2}
- \lambda_1 \delta(\delta+1)x^{-\delta-2} \cr
 \ar=\ar
\lambda_1\delta x^{-2\delta-2}[\lambda_1\delta
- (\delta+1)x^{\delta}] \cr
 \ar\ge\ar
\lambda_1\delta x^{-2\delta-2}[\lambda_1\delta
- (\delta+1)c_2^{\delta}] \cr
 \ar\ge\ar
\lambda_1^2\delta x^{-2\delta-2},
\qquad 0<x<c_2,\,y>0
 \eeqlb
and
 \beqlb\label{4.6}
g''_{xx}(x,y)/g(x,y)
 =
2\lambda_1(x-c^*)^{-6}[2\lambda_1
- 3(x-c^*)^2]
>0,
\quad c_3<x<c^*,\,y>0.
 \eeqlb
In addition, since   $g_0,g_0'$ and $g_0''$ are bounded on $[c_1,c_3]$,
then
 \beqlb\label{4.8}
C_0:=\sup_{x\ge c_1,y>0}\big[g(x,y)+|g_x'(x,y)|+|g_{xx}''(x,y)|\big]<\infty.
 \eeqlb

{\bf Step 2.}
In this step, we estimate $L_1g(x,y)$ which is defined in \eqref{3.21}.
Recall \eqref{3.1}.
Observe that $g(u,y)\ge g(x,y)$ for all $0<x\le u<c_2$ and $y>0$.
It follows from \eqref{3.1}, \eqref{7.2} and \eqref{4.5} that
for $x<c_1$, $z< c_2-c_1$ and $0<y<c^*$ we have
 \beqnn
 \ar\ar
g(x,y)^{-1}K_z^1g(x,y)
=z^2g(x,y)^{-1}\int_0^1g_{xx}''(x+zu,y)(1-u)\dd u \cr
 \ar\ar\qquad
\ge
\lambda_1^2z^2g(x,y)^{-1}
\int_0^1(x+zu)^{-2\delta-2}g(x+zu,y)(1-u)\dd u \cr
 \ar\ar\qquad
\ge
\lambda_1^2 z^2  \int_0^1(x+zu)^{-2\delta-2}(1-u)\dd u \cr
 \ar\ar\qquad
\ge
\lambda_1^2 x^{-2\delta-2}z^2 \int_0^x(1+c_2-c_1)^{-2\delta-2}(1-u)\dd u
\ge
\lambda_1^2 C_1 x^{-2\delta-1}z^2
 \eeqnn
for some constant $C_1>0$ independent of $\lambda_1$ and $\lambda_2$,
which gives
 \beqlb\label{5.6}
\int_0^{c_2-c_1}
K_z^1g(x,y)\mu(\dd z)
\ge
\lambda_1^2C_1\int_0^{c_2-c_1}z^2\mu(\dd z) x^{-2\delta-1}g(x,y),\quad
x\le c_1,~y>0.
 \eeqlb

Since $\lambda_1,\delta>1$, then by \eqref{3.1} and \eqref{6.5},
 \beqnn
K_z^1g(x,y)
 \ar\ge\ar
-g(x,y)-zg'_x(x,y)
= -g(x,y)-\lambda_1\delta z x^{-\delta-1}g(x,y) \cr
 \ar\ge\ar
-\lambda_1\delta x^{-\delta-1}g(x,y)(1+z),
\qquad x\le c_1,~y>0,
 \eeqnn
which implies that
 \beqlb\label{5.7}
\int_{c_2-c_1}^\infty
K_z^1g(x,y)\mu(\dd z)
\ge
-\lambda_1\delta x^{-\delta-1}g(x,y)
\int_{c_2-c_1}^\infty (1+z)\mu(\dd z),\quad
x\le c_1,~y>0.
 \eeqlb
Combining \eqref{5.6} and \eqref{5.7} we get
 \beqnn
g(x,y)^{-1}\int_0^\infty K_z^1g(x,y)\mu(\dd z)
\ge
\lambda_1^2 C_1\int_0^{c_2-c_1}z^2\mu(\dd z) x^{-2\delta-1}
-\lambda_1\delta x^{-\delta-1}\int_{c_2-c_1}^\infty (1+z)\mu(\dd z)
 \eeqnn
for $0<x<c_1$ and $0<y<c^*$.
Therefore, by \eqref{6.5} and \eqref{4.5},
 \beqlb\label{3.9a}
 \ar\ar
g(x,y)^{-1}L_1g(x,y) \cr
 \ar\ge\ar
\lambda_1 x^{-\delta}\Big[- \delta a_1(x)x^{-1}+\lambda_12^{-1}\delta a_2(x)x^{-2-\delta} \cr
 \ar\ar
+\Big(\lambda_1C_1\int_0^{c_2-c_1}z^2\mu(\dd z)
-\delta x^\delta \int_{c_2-c_1}^\infty (1+z) \mu(\dd z)\Big)
a_3(x)x^{-\delta-1}\Big]
 \eeqlb
for all $0<x< c_1, 0<y<c^*$.
Under condition (C2), $C_2:=\sup_{0<x<c^*}a_1(x)x^{-1}<\infty$.
Combining this with \eqref{3.9a} and \eqref{6.10}, for all $0<x< c_1, 0<y<c^*$
and large enough $\lambda_1$ with
 \beqnn
2^{-1}\lambda_1C_1\int_0^{c_2-c_1}z^2\mu(\dd z)
>
\delta c_1^\delta \int_{c_2-c_1}^\infty (1+z) \mu(\dd z),
 \eeqnn
we get
 \beqlb\label{3.9aa}
 \ar\ar
g(x,y)^{-1}L_1g(x,y) \cr
 \ar\ge\ar
\lambda_1 c_1^{-\delta}\Big[- C_2\delta+\lambda_12^{-1}\delta a_2(x)x^{-2-\delta} \cr
 \ar\ar\qquad\quad
+\Big(\lambda_1C_1\int_0^{c_2-c_1}z^2\mu(\dd z)
-\delta c_1^\delta \int_{c_2-c_1}^\infty (1+z) \mu(\dd z)\Big)
a_3(x)x^{-\delta-1}\Big] \cr
 \ar\ge\ar
\lambda_1 c_1^{-\delta}\Big[- C_2\delta+2^{-1}\lambda_1a_2(x)x^{-2-\delta}
+2^{-1}\lambda_1C_1\int_0^{c_2-c_1}z^2\mu(\dd z)
a_3(x)x^{-\delta-1}\Big] \cr
 \ar\ge\ar
\lambda_1 c_1^{-\delta}\Big[- C_2\delta+2^{-1}\lambda_1
\Big[1\wedge\Big(C_1\int_0^{c_2-c_1}z^2\mu(\dd z)\Big)\Big]
\big[a_2(x)x^{-2-\delta}
+a_3(x)x^{-\delta-1}\big]\Big] \cr
 \ar\ge\ar
\lambda_1 c_1^{-\delta}\Big[- C_2\delta+2^{-1}\lambda_1 c_0
\Big[1\wedge\Big(C_1\int_0^{c_2-c_1}z^2\mu(\dd z)\Big)\Big]\Big].
 \eeqlb
Observe that the term in the bracket of the above inequality is positive for large enough $\la_1$.
Thus for all large enough $\la_1>0$ there is a constant $d_1:=d_1(\lambda_1)>0$ so that
 \beqlb\label{3.9}
g(x,y)^{-1}L_1g(x,y)\ge d_1,\qquad 0<x< c_1,~0<y<c^*.
 \eeqlb
Since $g(x,y)=0$ for all $x\ge c^*$ or $y\ge c^*$, then $L_1g(x,y)=0$ for all $x\ge c^*$ or $y\ge c^*$.
By \eqref{6.5} and \eqref{4.6},
for large enough $\la_1$,
 \beqnn
-g'_x(x,y)=2\lambda_1(c^*-x)^{-3} g(x,y),~
g''_{xx}(x,y)
\ge
2\lambda_1^2(x-c^*)^{-6}g(x,y),\quad
c_3<x<c^*,~y>0.
 \eeqnn
Thus, for large enough $\la_1>0$ we have $-g'_x(x,y)\ge0$, $g''_{xx}(x,y)\ge0$, and
 \beqnn
K_z^1g(x,y)
=z^2\int_0^1g_{xx}''(x+zu,y)(1-u)\dd u\ge0
 \eeqnn
for all $x>c_3$ and $y>0$.
Now by the definition of $L_1g(x,y)$ in \eqref{3.21},
 \beqlb\label{3.9bb}
L_1g(x,y)\ge0,\qquad x\ge c_3,\,y>0
 \eeqlb
for large enough $\la_1>0$.
By \eqref{4.8}, for each $c_1\le x\le c_3$, $0<y<c^*$ and $z>0$,
 \beqnn
|K_z^1g(x,y)|
=z^2\Big|\int_0^1g_{xx}''(x+zu,y)(1-u)\dd u\Big|
\le C_0z^2
 \eeqnn
and
 \beqnn
K_z^1g(x,y)\ge-g(x,y)-zg'(x,y)\ge -C_0(1+z).
 \eeqnn
Then by \eqref{4.8} and the definition of $L_1g(x,y)$ in \eqref{3.21} again, for $c_1\le x\le c_3$ and $0<y\le c^*$,
 \beqlb\label{3.19}
L_1g(x,y)\ge -C_0 \Big[a_1(x)+2^{-1}a_2(x)+a_3(x)\int_0^1z^2\mu(\dd z)
+a_3(x)\int_1^\infty (1+z)\mu(\dd z)\Big].
 \eeqlb
Since $\inf_{c_1\le x\le c_3,\,0<y<c^*}g(x,y)>0$,
then by \eqref{3.19} and condition (C1),
it is elementary to see that
there is a constant $d_2:=d_2(\lambda_1)>0$ so that
 \beqlb\label{3.199}
g(x,y)^{-1}L_1g(x,y)\ge -d_2,\qquad c_1\le x\le c_3,\,0<y<c^*.
 \eeqlb
Combining the above inequality with \eqref{3.9} and \eqref{3.9bb}
we obtain
 \begin{equation}\label{3.9b}
 \begin{split}
&L_1g(x,y)\ge -d_2 g(x,y),\qquad x\ge c_1,\,0<y<c^*,\\
&L_1g(x,y)\ge d_1 g(x,y),\qquad 0<x< c_1,~0<y<c^*.
 \end{split}
 \end{equation}

{\bf Step 3.}
In this step we first estimate $L_2g(x,y)$ defined in \eqref{3.22} and then finish the proof.
By \eqref{6.5} and \eqref{4.7},
 \beqlb\label{3.24}
-g'_y(x,y)\ge0,~g''_{yy}(x,y)\ge0,\quad
0< x\le c^*,\,0<y<c^*,\,\lambda_2\ge 4r^{-1}.
 \eeqlb
Moreover, by \eqref{3.1} and \eqref{7.2},
 \beqlb\label{3.25}
K_z^2g(x,y)=z^2\int_0^1g''_{yy}(x,y+zv)(1-v)\dd v\ge0,\quad
0< x\le c^*,\,0<y<c^*,\,\lambda_2\ge 4r^{-1}.
 \eeqlb
Therefore, by the assumption of $0<r<1-\theta$, \eqref{6.9} and \eqref{6.5} again,
 \beqlb\label{3.23}
L_2g(x,y)
 \ar\ge\ar
-\kappa(x)\theta(y)g'_y(x,y)
=\lambda_2r\bar{c}\kappa(x)\theta(y)(\sin \bar{c}y)^{r-1} (\cos \bar{c}y)^{-1-r}
g(x,y) \cr
 \ar\ge\ar
\lambda_2r\bar{c}\kappa(x)\theta(y)(\bar{c}y)^{r-1}
g(x,y)
\ge
\lambda_2r\bar{c}^rc_0^2y^{\theta-1+r} g(x,y)
\ge 2d_2g(x,y)
 \eeqlb
for all $c_1\le x\le c^*$, $0<y<c^*$ and for $\lambda_2$ large enough,
where the constant $d_2>0$ is determined in \eqref{3.199}.
Since $g(x,y)=0$ for $x\ge c^*$ or $y\ge c^*$, then $L_2g(x,y)=0$ for $x\ge c^*$  or $y\ge c^*$.
It follows from \eqref{3.23} that
 \beqlb\label{3.26}
L_2g(x,y)
\ge 2d_2g(x,y),
\quad x\ge c_1,~y>0
 \eeqlb
and $L_2g(x,y)\ge0$ for all $x,y>0$ by \eqref{3.24} and \eqref{3.25}.
Recalling \eqref{3.20}.
Combining \eqref{3.26} with \eqref{3.9b} we get
 \beqnn
Lg(x,y)
 \ar=\ar
L_1g(x,y)+L_2g(x,y)\ge d_1 g(x,y),
\quad 0<x<c_1,\, y>0, \cr
Lg(x,y)
 \ar\ge\ar
[2d_2-d_2] g(x,y)
= d_2 g(x,y),
\qquad~~~ x\ge c_1,\, y>0,
 \eeqnn
which verifies  condition (iii) of Proposition \ref{t0.3}.

Therefore, by Proposition \ref{t0.3},
$\mbf{P}\{\tau_0^X\wedge\tau_0^Y<\infty\}\ge g(x_0,y_0)/[\sup_{x,y>0}g(x,y)]$ for large enough $\lambda_1,\lambda_2>0$ and $X_0,Y_0\in(0,c^*)$.
Since $\tau_0^X=\infty$ almost surely by Lemma \ref{t2.1} (ii),
we have $\mbf{P}\{\tau_0^Y<\infty\}>0$
for $0<X_0,Y_0<c^*$.
For general initial values $X_0>c^*$ or $Y_0>c^*$,
let $\tau:=\tau_{c^*}^{X+Y}$.
By Lemma \ref{t2.1} we have $\tau<\infty$ almost surely
and then by the Markov property,
 \beqnn
\mbf{P}\{\tau_0^Y<\infty\}
=\mbf{P}\{\tau_0^Y<\infty|(X_\tau,Y_\tau)\}>0,
 \eeqnn
which completes the proof.
\qed

\subsection{Proof of Theorem \ref{t0.01}}

\blemma\label{t1.4}
Suppose that Condition \ref{c} (ib)-(ic) hold.
Let $\tilde{g}$ be a nonnegative process
satisfying $\int_0^\infty \tilde{g}(s)^{\delta}\dd s=\infty$ almost surely for some constant
$\delta$ with
$q/(q+1-\theta)<\delta\le1$.
Let $(u_t)_{t\ge0}$ be the non-negative solution to
 \beqnn
u_t
 \ar=\ar
u_0
-\int_0^t[b_1(u_s)+\theta(u_s)\tilde{g}(s)]\dd s \cr
 \ar\ar
+\int_0^tb_2(u_s)^{1/2}\dd W_s
+\int_0^t\int_0^\infty\int_0^{b_3(u_{s-})}z \tilde{N}(\dd s,\dd z,\dd u).
 \eeqnn
Then for each $u_0>0$, we have
$\mbf{P}\{\tau_0^u<\infty\}=1$.
\elemma
\proof
To establish the proof we apply Corollary \ref{cor2.3}.
The key is to construct the function $g$ and verify the conditions
in Corollary \ref{cor2.3}.
For $r\in(0,1-\theta)$ and $v,\lambda>0$ let $g(v)=\e^{-\lambda v^r}$.
Then
 \beqlb\label{6.7}
g'(v)=-r\lambda v^{r-1}g(v),~
g''(v)
=
r\lambda[r\lambda v^r+(1-r)]v^{r-2}g(v)
\ge r(1-r)\lambda v^{r-2}g(v).
 \eeqlb
By It\^o's formula we can see that the operator $L_t$ is given by
 \beqnn
L_tg(v)
:=
-[b_1(v)+\theta(v)\tilde{g}(t)]g'(v)
+
2^{-1} b_2(v)g''(v)
+b_3(v)\int_0^\infty K_zg(v)\nu(\dd z),
 \eeqnn
where $K_zg(v)$ is given in \eqref{7.1}.
In the following we find an estimation of $L_tg(v)$.
It follows from \eqref{7.7b} that
 \beqlb\label{6.1}
\int_0^\infty  K_zg(v)\nu(\dd z)
\ge
\lambda r(1-r)v^rH_{2,0}(v)g(v),
 \eeqlb
where the function $H_{2,0}$ is defined in \eqref{0.6}.
Therefore, by \eqref{6.7} and \eqref{6.1}, for all $n\ge1$,
 \beqlb\label{4.1}
L_tg(v)
 \ar\ge\ar
\la rg(v)[b_1(v)v^{r-1}
+(1-r)2^{-1} b_2(v)v^{r-2}
+(1-r)b_3(v)v^rH_{2,0}(v)] \cr
 \ar\ge\ar
\la r(1-r)g(v)v^rG_{2,0}(v) \ge
\la r(1-r)|c^*|^rd_ng(v),\qquad c^*\le v<n,
 \eeqlb
where the function $G_{2,0}$ is defined in \eqref{0.8} and
$d_n:=\inf_{c^*\le v<n}G_{2,0}(v)>0$.
Under Condition \ref{c} (ib)-(ic), we have
 \beqnn
L_tg(v)
 \ar\ge\ar
\la r v^rg(v)\big[b_1(v)v^{-1}
+(1-r)2^{-1} b_2(v)v^{-2}
+(1-r)b_3(v)H_{2,0}(v)
+c_\theta\tilde{g}(t)v^{\theta-1}\big] \cr
 \ar\ge\ar
\la r(1-r) v^rg(v)\big[G_{2,0}(v)
+c_\theta\tilde{g}(t)v^{\theta-1}\big] \cr
 \ar\ge\ar
\la r(1-r) g(v)\big[bv^{q+r}+c_\theta\tilde{g}(t)v^{\theta-1+r}\big],
\qquad 0<v\le c^*.
 \eeqnn
Then
 \beqlb\label{4.2''}
L_tg(v)
\ge
\la r(1-r)c_\theta g(v)\tilde{g}(t)|c^*|^{\theta-1+r},
\qquad 0<v\le c^*,
 \eeqlb
and by Lemma \ref{t2.3}, there are constants $C_1=C_1(r)>0$
and $C_2=C_2(r)>0$ so that
 \beqlb\label{4.2'}
g(v)^{-1}L_tg(v)
 \ar\ge\ar
C_1\la v^{(1-1/\bar{q})(r+q)+(\theta-1+r)/\bar{q}}\tilde{g}(t)^{1/\bar{q}} \cr
 \ar=\ar
C_1\la \tilde{g}(t)^{1/\bar{q}}v^{r+q-(q+1-\theta)/\bar{q}}
\ge C_2\la \tilde{g}(t)^{1/\bar{q}},
\quad 0<v\le c^*
 \eeqlb
for $\bar{q}>1$ and $r+q-(q+1-\theta)/\bar{q}\le0$, which
is equivalent to
 \beqnn
\frac{1}{\bar{q}}\ge\frac{r+q}{q+1-\theta}.
 \eeqnn
It holds as long as
$r$ is small enough and
 \beqnn
\frac{1}{\bar{q}}>\frac{q}{q+1-\theta}.
 \eeqnn

Combining Corollary \ref{cor2.3} and
\eqref{4.1}-\eqref{4.2'} one gets
$\mbf{P}\{\tau_0^u<\infty\}\ge \e^{-\lambda u_0^r}$
if either
 \beqnn
\int_0^\infty \tilde{g}(s)\dd s=\infty, \qquad
\mbox{or }\int_0^\infty \tilde{g}(s)^{1/\bar{q}}\dd s=\infty~~
\mbox{for }1>\frac{1}{\bar{q}}>\frac{r+q}{q+1-\theta}.
 \eeqnn
Taking $\bar{q}=1/\delta$ and
letting $\lambda\to0$ we get $\mbf{P}\{\tau_0^u<\infty\}=1$ under the above conditions.
This finishes the proof.
\qed

\blemma\label{t1.5}
Suppose that Condition \ref{c}(ia) holds.
Then for $0<\bar{p}<1$ and $\kappa>0$ satisfying $\kappa \bar{p}\le p$, we have
$\int_0^\infty X_s^{\kappa\bar{p}}\dd s=\infty$
almost surely.
\elemma
\proof
Let $\bar{a}_i(x)=a_i(x)/x^{\kappa \bar{p}}$ for $i=1,2,3$.
Then by the same argument as in \cite[Theorem 2.15]{LYZh},
there are, on an extended  probability space, a Brownian motion
$(\bar{B}_t)_{t\ge0}$ and compensated Poisson random measure
$\{\tilde{\bar{M}}(\dd t,\dd z , \dd u)\}$ with intensity
%are $\dd t\mu(\dd z)\dd u$ so that there is a nonnegative process
 $\dd t\mu(\dd z)\dd u$ so that there is a nonnegative process
$(\bar{X}_t)_{t\ge0}$ solving:
 \beqnn
\bar{X}_t=\bar{X}_0 -\int_0^t\bar{a}_1(\bar{X}_s)\dd s
+\int_0^t\bar{a}_2(\bar{X}_s)^{1/2}\dd \bar{B}_s
+\int_0^t\int_0^\infty
\int_0^{\bar{a}_3(\bar{X}_{s-})}z \tilde{\bar{M}}(\dd s,\dd z,\dd u).
 \eeqnn
Moreover,
by \cite[Proposition 2.16]{LYZh},
 \beqlb\label{3.6}
\int_0^\infty X_s^{\kappa \bar{p}}\dd s=\tau_0^{\bar{X}}
 \eeqlb
almost surely.
Recall the function $H_{1,0}$ in \eqref{0.5}.
Since $\kappa \bar{p}\le p$, then under condition \ref{c}(ia),
for all $0<u<1$,
 \beqnn
\bar{a}_1(u)u^{-1}+2^{-1}\bar{a}_2(u)u^{-2}+\bar{a}_3(u)H_{1,0}(u)
=G_{1,0}(u)u^{-\kappa \bar{p}}\le G_{1,0}(u)u^{-p}
\le a.
 \eeqnn
Now by Lemma \ref{t2.1} (ii),
$\tau_0^{\bar{X}}=\infty$ almost surely.
Then the assertion follows from \eqref{3.6} immediately.
\qed

\noindent{\it Proof of Theorem \ref{t0.01}.}
Let $\bar{p}\in(0,1)$ satisfy $\bar{p}\kappa\le p$.
By Lemmas \ref{t1.5} and \ref{t2.1}(i) and Condition \ref{c} (ic),
 \beqnn
\int_0^\infty \kappa(X_s)^{\bar{p}}\dd s=\infty
 \eeqnn
almost surely.
Taking $\delta=\bar{p}$ in Lemma \ref{t1.4}, we have
$\mbf{P}\{\tau_0^Y<\infty \}=1$
for
 \beqnn
\frac{q}{q+1-\theta}<\bar{p} ~~
\mbox{and}~~\bar{p}\kappa\le p,
 \eeqnn
which finishes the proof.
\qed

\subsection{Proof of Theorem \ref{t0.02}}

Different from proofs of the previous theorems, in Theorem \ref{t0.02}   we adopt an approach similar to that in \cite{LYZh}.   In the proof  we first identify a power function $\hat{X}_t$ such that on one hand, with a positive probability $\hat{X}_t$ stays above $X_t$ for all large $t$, and on the other hand,  $\hat{X}_t$ decreases faster enough for large $t$ so that if the term $\kappa(X_s)$  is replaced by $\kappa(\hat{X}_t)$  in the SDE for $Y$ in  \eqref{1.1}, then $Y$ becomes extinguishing with a strictly positive probability. To implement this idea we further approximate $\hat{X}_t$ by a step function and construct a process $\hat{Y}_t$ as a piecewise solution to a modified SDE for $Y$.
The desired result then follows from a comparison theorem for SDE.

We state the following comparison theorem with its proof postponed to the Appendix.
For $i=1,2$ let $\{B_i(t,u):t\ge0,u\in\mbb{R}\}$
be a two-parameter real-valued process
with  $(u,\omega)\mapsto B_i(t,u,\omega)$
measurable with respect to
$\mathscr{B}(\mbb{R})\times\mathscr{F}_t$ for each $t\ge0$.
Let $U$ and $V$ be Borel functions on $\mbb{R}$
and $V\ge0$.
\bproposition\label{t4.4}
For $i=1,2$, let the c\`adl\`ag $\mbb{R}$-valued process $(x_i(t))_{t\ge0}$ be the solution to SDE
 \beqlb\label{8.1}
x_i(t)
 \ar=\ar
x_i(0)
+\int_0^tB_i(s,x_i(s))\dd s
+\int_0^t U(x_i(s))\dd W_s \cr
 \ar\ar
+\int_0^t\int_0^\infty\int_0^{V(x_i(s-))}z \tilde{N}(\dd s,\dd z,\dd u).
 \eeqlb
Suppose that
$B_1(t,u)\le B_2(t,u)$ for all $t\ge0$ and $u\in\mbb{R}$.
In addition, assume that $u\mapsto V(u)$ is nondecreasing
and that
there exists a sequence of increasing stopping times $(\gamma_n)_{n\geq 1}$ and a sequence of nonnegative constants $(C_n)_{n\geq 1}$ so that
 \beqnn
|B_1(s,u)-B_1(s,v)|+|U(u)-U(v)|+|V(u)-V(v)|\le C_n|u-v|
 \eeqnn
for all $n^{-1}\le|u|,|v|\le n$ and $s\le \gamma_n$.
If $\mbf{P}\{x_1(0)\le x_2(0)\}=1$,
then
\[\mbf{P}\{x_1(t)\le x_2(t)\quad\text{ for all}\quad 0<t<\tilde{\gamma}\}=1,\]
where
 \beqlb\label{4.2}
 \tilde{\gamma}:=\lim_{n\to\infty}\tilde{\gamma}_n
\quad\text{and}\quad
\tilde{\gamma}_n:=\gamma_n\wedge\tau_{1/n}^{x_1}\wedge\tau_{1/n}^{x_2}
\wedge\sigma_n^{x_1}\wedge\sigma_n^{x_2}
 \eeqlb
with $\tau_{1/n}^{x_i}:=\inf\{t\ge0: |x_i(t)|\le 1/n\}$
and $\sigma_n^{x_i}:=\inf\{t\ge0: |x_i(t)|\ge n\}$
for $i=1,2$.
\eproposition

Recall the constant $c^*$ in Condition \ref{c} and the definitions of stopping times
$\tau_w^X=\tau^X(w)$ and $\sigma_w^X=\sigma^X(w)$ in \eqref{1.2} and \eqref{1.3}
for constant $w>0$.

\blemma\label{t4.1}
Under Condition \ref{c} (iia) with $p>0$,
for any $0<w<X_0$ and $0<v\le c^*$ we have
 \beqnn
\mbf{E}\big[\tau_w^X\wedge \sigma_v^X\big]
\le
2p(p\wedge1)^{-1}[1-2^{-1}(p\wedge1)]^{-1}a^{-1}(X_0-w) w^{-p-1}.
 \eeqnn
\elemma
\proof
It is elementary to see that for $\delta>\delta_1>0$ and $u>0$,
 \beqnn
(1+u)^{-\delta}-1+\delta u
\ge
(1+u)^{-\delta_1}-1+\delta_1 u
\ge
-[(1+u)^{\delta_1}-1-\delta_1 u]
\ge0,
 \eeqnn
which implies that for $p_1:=2^{-1}(p\wedge1)$,
 \beqlb\label{5.0}
p(p+1)H_{1,p}(u)
\ge p_1(1-p_1)H_{1,-p_1}(u)
\ge p_1(1-p_1)H_{1,0}(u),
 \eeqlb
where the function $H_{1,p}$ and $H_{1,0}$ are defined in \eqref{0.3}
and \eqref{0.5}, respectively.

By \eqref{1.1} and It\^o's formula,
 \beqnn
X_t^{-p}
 \ar=\ar
X_0^{-p}+p\int_0^ta_1(X_s)X_s^{-p-1}\dd s+\frac{p(p+1)}{2} \int_0^ta_2(X_s)X_s^{-p-2}\dd s\cr
 \ar\ar
+p(p+1)\int_0^ta_3(X_s)X_s^{-p}H_{1,p}(X_s)\dd s
-p\int_0^t a_2(X_s)^{1/2}X_s^{-p-1}\dd B_s \cr
 \ar\ar
+\int_0^t\int_0^\infty\int_0^{a_3(X_{s-})}
[(X_{s-}+z )^{-p}-X_{s-}^{-p}]\tilde{M}(\dd s,\dd z, \dd u).
 \eeqnn
Since
$X_t^{-p}\le w^{-p}$ for $0\le t\le \tau_w^X$,
then by \eqref{5.0}, for $0<p_1<p\wedge1$, we have
 \beqnn
w^{-p}
 \ar\ge\ar
\mbf{E}\Big[X_{t\wedge \tau_w^X\wedge \sigma_v^X}^{-p}\Big] \cr
\ar\ge\ar
X_0^{-p}+\mbf{E}\Big[\int_0^{t\wedge \tau_w^X\wedge \sigma_v^X
\Big(pa_1(X_s)X_s^{-p-1}+\frac{p(p+1)}{2} a_2(X_s)X_s^{-p-2}} \cr
 \ar\ar\qquad\qquad\qquad\qquad\qquad\quad
+p(p+1)a_3(X_s)X_s^{-p}H_{1,p}(X_s)\Big)\dd s\Big] \cr
 \ar\ge\ar
X_0^{-p}+p_1(1-p_1)\mbf{E}\Big[\int_0^{t\wedge \tau_w^X\wedge \sigma_v^X}
G_{1,0}(X_s)X_s^{-p}\dd s\Big] \cr
 \ar\ge\ar
X_0^{-p}
+p_1(1-p_1)a\mbf{E}\big[t\wedge \tau_w^X\wedge \sigma_v^X\big],
 \eeqnn
where Condition \ref{c} (iia) is used in the last inequality.
Using Fatou's lemma we get
 \beqnn
\mbf{E}[\tau_w^X\wedge \sigma_v^X]
 \ar\le\ar
\liminf_{t\to\infty}\mbf{E}\big[t\wedge \tau_w^X\wedge \sigma_v^X\big] \cr
 \ar\le\ar
(p_1(1-p_1)a)^{-1}(w^{-p}-X_0^{-p}) \cr
 \ar\le\ar
p(p_1(1-p_1)a)^{-1}(X_0-w)w^{-p-1},
 \eeqnn
where we need the mean value theorem for  the last inequality.
This ends the proof.
\qed

Since process $X$ turns to decrease in the long run,
we next find a power function of time that is uniformly larger than $X_t$ for all large $t$ with a probability close to one. To this end, we consider a partition of the duration of time into consecutive time intervals  with  partition points increasing geometrically. Then for an arbitrary time interval in the partition, we further choose three levels $0< l_1<l_2<l_3$ properly so  that, during this time interval, process $X$ typically reaches level $l_1$ first before upcrossing level $l_3$, and then it typically stays below level $l_2$ continuously after having reached level $l_1$. The above choices of time partition and the associated levels allow us to show that process $(X_t)_{t\geq 0}$ typically stay below the desired power function of $t$ for all large $t$.

\blemma\label{t4.2}
Under Condition \ref{c} (iia) with $p>0$, for any $\delta>0$ and small enough $\varepsilon\in(0,1)$,
there are constants $C(\delta,\varepsilon)>0$ and $\delta_1\in(0,1)$ that does not depend on  $\varepsilon$
so that for $X_0=\varepsilon^m$ with large enough $m\ge1$, we have
 \beqnn
\mbf{P}\{X_t\le t^{-\frac{1}{p+\delta}}\wedge c^* \,\,\, \text{for all}\,\,\, t>0 \}
\ge
1-C(\delta,\varepsilon)\varepsilon^{m\delta_1/8}.
 \eeqnn
\elemma
\proof
In the following let $\varepsilon_n:=\varepsilon^n$ for $n\geq1$.
For any $\delta>0$, let
\beqlb\label{5.5}
\delta_1:=\frac{\delta}{2p+2+\delta}<1.
 \eeqlb
For any fixed positive integer $m$ define
 \beqnn
\bar{K}_m:=\Big\{\sup_{ t\leq \varepsilon_m^{-p -(p +2)\delta_1}}X_t\leq \varepsilon_m^{1-\delta_1}, \, X_{\varepsilon_m^{-p -(p +2)\delta_1}}\leq
\varepsilon_{m+1} \Big\}
 \eeqnn
and
 \beqnn
\bar{K}_n:=\Big\{\sup_{\varepsilon_{n-1}^{-p -(p +2)\delta_1}\leq t< \varepsilon_n^{-p -(p +2)\delta_1}}X_t\leq \varepsilon_n^{1-\delta_1}, \,
X_{\varepsilon_n^{-p -(p +2)\delta_1}}\leq \varepsilon_{n+1} \Big\}
 \eeqnn
for $n>m$. In the following we first show that $X_t\le t^{-1/(p +\delta)}\wedge c^*$
for all $t>0$ on $\cap_{n=m}^\infty \bar{K}_n$.

It is obvious that $X_t\le c^*$ for all $t\ge0$
on event $\cap_{n=m}^\infty \bar{K}_n$ for  $\varepsilon$ small enough.
Let $r=\frac{1}{p+\delta}$.
By \eqref{5.5}, it is easy to check that
 \beqlb\label{5.10}
-\frac{1}{p+\delta}=-r=-\frac{1-\delta_1}{p +(p +2)\delta_1}.
 \eeqlb
Thus, for all $n\geq m$ we have
$(\varepsilon_n^{-p -(p +2)\delta_1})^{-r}=\varepsilon_n^{1-\delta_1}$
by \eqref{5.10}.
Therefore, on event $\cap_{n=m}^\infty \bar{K}_n$, for any $t>0$ with
 \beqlb\label{5.11}
\varepsilon_{n-1}^{-p-(p+2)\delta_1}\leq t<  \varepsilon_n^{-p-(p +2)\delta_1},\quad
n\geq m+1,
 \eeqlb
we have
 \beqnn
X_t\leq \varepsilon_n^{1-\delta_1}=\varepsilon_{n}^{-r(-p -(p +2)\delta_1)}< t^{-r}
 \eeqnn
by \eqref{5.10} again and \eqref{5.11},
and for $0\leq t\leq  \varepsilon_m^{-p -(p +2)\delta_1}$ we also have
 \beqnn
X_t\leq \varepsilon_m^{1-\delta_1}= \varepsilon_{m}^{-r(-p-(p+2)\delta_1)}\le t^{-r}.
 \eeqnn

We now  estimate the probability of $\cap_{n=m}^\infty \bar{K}_n$.
In the rest of the proof we use notations
 \beqnn
\mbf{E}_{\bar{\varepsilon}}[\,\cdot\,]
=\mbf{E}\big[\cdot|X_0=\bar{\varepsilon}\big]
~\mbox{ and }~
\mbf{P}_{\bar{\varepsilon}}\{\,\cdot\,\}
=\mbf{P}\big\{\cdot|X_0=\bar{\varepsilon}\big\},
\qquad \bar{\varepsilon}>0.
 \eeqnn
By Lemma \ref{t4.1},
there is a constant $c_1>0$ independent of $\varepsilon$ and $n$
so that
 \beqlb\label{5.12}
\mathbf{E}_{\bar{y}}
\Big[\tau^X(\varepsilon_{n+1}^{1+\delta_1})\wedge\sigma_{c^*}^X\Big]
\le
c_1a^{-1}(\varepsilon_n-{\varepsilon_{n+1}^{1+\delta_1}})
\varepsilon_{n+1}^{(1+\delta_1)(-p -1)}
\le
c_1a^{-1}\varepsilon^{-1}\varepsilon_{n+1}^{-p-(p +1)\delta_1}
 \eeqlb
for all $0<\bar{y}\le \varepsilon_n$,
where the constant $a>0$ is determined in Condition \ref{c} (iia).
Using the Markov inequality and \eqref{5.12} we obtain
 \beqlb\label{5.1}
 \ar\ar
\mbf{P}_{\bar{y}}\Big\{
\tau^X(\varepsilon_{n+1}^{1+\delta_1})\wedge\sigma_{c^*}^X
>\varepsilon_n^{-p -(p +2)\delta_1}-\varepsilon_{n-1}^{-p -(p +2)\delta_1}\Big\} \cr
 \ar\ar\qquad\le
[\varepsilon_n^{-p -(p +2)\delta_1}-\varepsilon_{n-1}^{-p -(p +2)\delta_1}]^{-1}
\mbf{E}_{\bar{y}}
\big[\tau^X(\varepsilon_{n+1}^{1+\delta_1})\wedge\sigma_{c^*}^X\big] \cr
 \ar\ar\qquad\le
c_1a^{-1}(1-\varepsilon^{p+(p+2)\delta_1})^{-1}\varepsilon^{-p-1 +(n-p -1)\delta_1}
 \eeqlb
for all $0<\bar{y}\le \varepsilon_n$.
By Lemma \ref{t3.2},
there is a constant $C>0$ independent of $\varepsilon$ and $n$
so that
 \beqlb\label{5.2}
\mbf{P}_{\bar{y}}\Big\{\sup_{t\ge0} X_t\geq \varepsilon_n^{\delta_1-1}\Big\}
\le C(\varepsilon_n^{\delta_1-1} {\bar y})^{1/4}
\le C\varepsilon^{n\delta_1/4}
 \eeqlb
for all $0<\bar{y}\le \varepsilon_n$.
By Fatou's lemma and \eqref{1.1},
 \beqnn
\varepsilon_n\mathbf{P}_{\varepsilon_n^{1+\delta_1}}
\big\{\sigma^X(\varepsilon_n)<\infty\big\}
\le
\mathbf{E}_{\varepsilon_n^{1+\delta_1}}\big[X_{\sigma^X(\varepsilon_n)}\big]
\le\liminf_{t\to\infty}
\mathbf{E}_{\varepsilon_n^{1+\delta_1}}\big[X_{t\wedge\sigma^X(\varepsilon_n)}\big]
\le\varepsilon_n^{1+\delta_1},
 \eeqnn
which implies
 \beqlb\label{5.3}
\mathbf{P}_{\varepsilon_n^{1+\delta_1}}\big\{\sigma^X(\varepsilon_n)<\infty\big\}
\le\varepsilon_n^{\delta_1}.
 \eeqlb
Similarly,
 \beqlb\label{5.8}
\mathbf{P}_{{\bar y}}\big\{\sigma^X_{c^*}<\infty\big\}
\le {c^*}^{-1}{\bar y}
\le {c^*}^{-1}\varepsilon_n
 \eeqlb
for all $\bar{y}\le \varepsilon_n$.

To consider the complements of events $\bar{K}_m$ and $\bar{K}_n$,  in the following, for the fixed $m\ge1$, we introduce
 \beqnn
E_m
 \ar:=\ar
\Big\{\sup_{  t\leq \varepsilon_m^{-p -(p +2)\delta_1}}X_t\geq \varepsilon_m^{1-\delta_1} \Big\}\bigcup
\Big\{\tau^X(\varepsilon_{m+1}^{1+\delta_1})\wedge\sigma_{c^*}^X
>\varepsilon_m^{-p -(p +2)\delta_1}\Big\} \cr
 \ar\ar
\qquad
\bigcup\Big\{\tau^X(\varepsilon_{m+1}^{1+\delta_1})\wedge\sigma_{c^*}^X< \infty, \sigma^X(\varepsilon_{m+1})
\circ\vartheta(\tau^X(\varepsilon_{m+1}^{1+\delta_1})\wedge\sigma_{c^*}^X)\leq \varepsilon_m^{-p -(p +2)\delta_1}\Big\}
 \eeqnn
and for $n>m$,
 \beqnn
E_n
 \ar:=\ar
\Big\{\sup_{\varepsilon_{n-1}^{-p -(p +2)\delta_1}\leq  t\leq \varepsilon_n^{-p -(p +2)\delta_1}}X_t\geq \varepsilon_n^{1-\delta_1}
\Big\} \cr
 \ar\ar\quad
\bigcup\Big\{\tau^X(\varepsilon_{n+1}^{1+\delta_1})\wedge\sigma_{c^*}^X
>\varepsilon_n^{-p -(p +2)\delta_1}-\varepsilon_{n-1}^{-p -(p +2)\delta_1}\Big\} \cr
 \ar\ar\quad
\bigcup\Big\{\tau^X(\varepsilon_{n+1}^{1+\delta_1})\wedge\sigma_{c^*}^X< \infty, \sigma^X_{\varepsilon_{n+1}}
\circ\vartheta(\tau^X(\varepsilon_{n+1}^{1+\delta_1})\wedge\sigma_{c^*}^X)\leq \varepsilon_n^{-p -(p +2)\delta_1}\Big\},
 \eeqnn
where $\vartheta(t)$ denotes the usual shift operator.
For $0<\bar{y}\leq \varepsilon_n$,  we have
 \beqnn
 \ar\ar
\mbf{P}_{\bar{y}}\Big\{\tau^X(\varepsilon_{n+1}^{1+\delta_1})\wedge\sigma_{c^*}^X< \infty, \sigma^X(\varepsilon_{n+1})
\circ\vartheta(\tau^X(\varepsilon_{n+1}^{1+\delta_1})\wedge\sigma_{c^*}^X)\leq \varepsilon_n^{-p -(p +2)\delta_1}\Big\} \cr
 \ar\le\ar
\mbf{P}_{\bar{y}}\{\sigma_{c^*}^X<\infty\}
+\mbf{P}_{\bar{y}}\Big\{\tau^X(\varepsilon_{n+1}^{1+\delta_1})< \infty, \sigma^X(\varepsilon_{n+1})
\circ\vartheta(\tau^X(\varepsilon_{n+1}^{1+\delta_1}))\leq \varepsilon_n^{-p -(p +2)\delta_1}\Big\} \cr
 \ar\le\ar
\mbf{P}_{\bar{y}}\{\sigma_{c^*}^X<\infty\}
+\mathbf{P}_{\varepsilon_{n+1}^{1+\delta_1}}\Big\{\sigma^X(\varepsilon_{n+1})\leq \varepsilon_n^{-p -(p +2)\delta_1}\Big\} \cr
 \ar\le\ar
\mbf{P}_{\bar{y}}\{\sigma_{c^*}^X<\infty\}
+\mathbf{P}_{\varepsilon_{n+1}^{1+\delta_1}}
\Big\{\sigma^X(\varepsilon_{n+1})<\infty\Big\},
 \eeqnn
and then by \eqref{5.1}-\eqref{5.8},
 \beqnn
 \ar\ar
\mbf{P}_{\bar{y}}\{E_n\} \cr
 \ar\leq\ar
\mbf{P}_{\bar{y}}
\Big\{\sup_{t\leq \varepsilon_n^{-p -(p +2)\delta_1}}X_t\geq \varepsilon_n^{1-\delta_1} \Big\}
+\mbf{P}_{\bar{y}}
\Big\{\tau^X(\varepsilon_{n+1}^{1+\delta_1})\wedge\sigma_{c^*}^X
>\varepsilon_n^{-p-(p+2)\delta_1}- \varepsilon_{n-1}^{-p-(p +2)\delta_1} \Big\}  \cr
 \ar\ar
+\mbf{P}_{\bar{y}}\Big\{\tau^X(\varepsilon_{n+1}^{1+\delta_1})\wedge\sigma_{c^*}^X< \infty, \sigma^X(\varepsilon_{n+1})
\circ\vartheta(\tau^X(\varepsilon_{n+1}^{1+\delta_1})\wedge\sigma_{c^*}^X)\leq \varepsilon_n^{-p -(p +2)\delta_1}\Big\}  \cr
 \ar\leq \ar
\mbf{P}_{\bar{y}}\Big\{\sup_{t\leq \varepsilon_n^{-p-(p+2)\delta_1}}X_t\geq \varepsilon_n^{1-\delta_1} \Big\}
+\mbf{P}_{\bar{y}}\{\sigma_{c^*}^X<\infty\}
+\mathbf{P}_{\varepsilon_{n+1}^{1+\delta_1}}
\Big\{\sigma^X(\varepsilon_{n+1})<\infty\Big\} \cr
 \ar\ar
+\mbf{P}_{\bar{y}}\Big\{\tau^X(\varepsilon_{n+1}^{1+\delta_1})
\wedge\sigma_{c^*}^X
>\varepsilon_n^{-p -(p +2)\delta_1}
-\varepsilon_{n-1}^{-p -(p +2)\delta_1}\Big\} \cr
 \ar\leq \ar
C\varepsilon_n^{\delta_1/4}
+\varepsilon_n/c^*+\varepsilon_{n+1}^{\delta_1}
+c_1a^{-1}(1-\varepsilon^{p+(p+2)\delta_1})^{-1}\varepsilon^{-p-1 +(n-p -1)\delta_1}
\cr
 \ar\leq \ar
\varepsilon^{n\delta_1/8}(1+\varepsilon^{-(p+1)(\delta_1+1)})=:M_n
 \eeqnn
for small enough $\varepsilon$.
Similarly, for small enough $\varepsilon$ and $0<\bar{y}\le \varepsilon_m$,
 \beqnn
\mbf{P}_{\bar{y}}\{E_m\}\leq \varepsilon^{m\delta_1/8}(1+\varepsilon^{-(p+1)(\delta_1+1)})=:M_m.
 \eeqnn
Let $\bar{K}_n^c$ denote the complement of set $\bar{K}_n$.
Note that $\bar{K}_n^c\subset E_n $.
Then $\mbf{P}_{\bar{y}}\{\bar{K}_n^c\}\le M_n$ for all $n\ge m$
and $0<\bar{y}\leq \varepsilon_n$.
It follows that
 \beqnn
\mbf{P}_{\varepsilon_m}\{\cup_{n=m}^\infty \bar{K}_n^c\}
 \ar=\ar
\mbf{P}_{\varepsilon_m}\{ \bar{K}_m^c\}+ \sum_{n=m+1}^\infty
\mbf{P}_{\varepsilon_m}\{\cap_{i=m}^{n-1}\bar{K}_i \cap \bar{K}_n^c\} \cr
 \ar=\ar
\mbf{P}_{\varepsilon_m}\{ \bar{K}_m^c\}+\sum_{n=m+1}^\infty \mbf{E}_{\varepsilon_m}
\Big[ 1_{\cap_{i=m}^{n-1}\bar{K}_i}
\mbf{P}\big\{ \bar{K}_n^c|X_{\varepsilon_{n-1}^{-p -(p +2)\delta_1}}\big\}
\Big] \cr
 \ar\le\ar
M_m+\sum_{n=m+1}^\infty M_n\mbf{E}_{\varepsilon_m}
\big[ 1_{\cap_{i=m}^{n-1}\bar{K}_i}
\big] \cr
 \ar\leq\ar
\sum_{n=m}^\infty M_n
= (1-\varepsilon^{\delta_1/8})^{-1}
(1+\varepsilon^{-(p+1 )(1+\delta_1)})\varepsilon^{m\delta_1/8}.
 \eeqnn
Then
 \beqnn
\mbf{P}_{\varepsilon_m}\{\cap_{n=m}^\infty \bar{K}_n\}= 1
-\mbf{P}_{\varepsilon_m}\{\cup_{n=m}^\infty \bar{K}_n^c\},
 \eeqnn
which finishes the proof.
\qed

Using the estimate as function of time obtained in  Lemma \ref{t4.2} we can construct a process $\hat{Y}$ which
does not become extinct with a positive probability and
can be shown by the comparison theorem  to be uniformly smaller than process $Y$
with a probability close to one.

For small enough $\delta,\epsilon\in(0,{c^*}\wedge1)$, suggested by Lemma \ref{t4.2}
we define
 \beqlb\label{7.12}
\hat{X}(t)= t^{-\frac{1}{p+\delta}}\wedge c^*
 \eeqlb
and $\epsilon_{n}=\epsilon^{n^2}$.
Let $(Y_1(t))_{t\ge0}$ be the nonnegative solution to
 \beqnn
Y_1(t)
 \ar=\ar
Y_0-\int_0^t\big[b_1(Y_1(s))+\theta(Y_1(s))\hat{X}(0)^\kappa\big]\dd s \cr
 \ar\ar
+\int_0^t b_2(Y_1(s))^{1/2}\dd W_s
+\int_0^t\int_0^\infty\int_0^{b_3(Y_1(s-))}z \tilde{N}(\dd s,\dd z,\dd u)
 \eeqnn
and $\gamma_1:=\inf\{t\geq 0: Y_1(t)<\epsilon_1\}$.
Define $\hat{Y}(t):=Y_1(t) $ for $t\in [0, \gamma_1]$.
Suppose that $\hat{Y}(t)$ has been defined for $t\in [0, T_n]$
with $T_n:=\sum_{i=1}^n \gamma_i$.
Let $(Y_{n+1}(t))_{t\ge0}$ be the nonnegative solution to
 \beqlb\label{7.3}
Y_{n+1}(t)
 \ar=\ar
Y_n(T_n)
-\int_0^t\big[b_1(Y_{n+1}(s))
+\theta(Y_{n+1}(s))\hat{X}(T_n)^\kappa\big]\dd s \cr
 \ar\ar
+\int_0^t b_2(Y_{n+1}(s))^{1/2}\dd W_s
+\int_0^t\int_0^\infty\int_0^{b_3(Y_{n+1}(s-))}z \tilde{N}(\dd s,\dd z,\dd u)
 \eeqlb
and $\gamma_{n+1}:=\inf\{t\geq 0: Y_{n+1}(t)<\epsilon_{n+1}\}$. Define $\hat{Y}(t):=Y_{n+1}(t-T_n)$ for $t\in (T_n, T_n+\gamma_{n+1}]=(T_n, T_{n+1}]$.
Then by the argument in \cite[Theorem 3.1]{LYZh}
and Condition \ref{c}(iii),
$\hat{Y}$ is a piecewise time homogeneous spectrally positive Markov process.

Choose $l$ satisfying that
 \beqlb\label{7.5}
0<l<q \quad\text{ and}\quad  \frac{l\kappa}{p +\delta}-l+\theta-1>0.
 \eeqlb
Such a value $l$ exists if \eqref{1.4} holds
and $\delta>0$ is small enough. In the next lemma,  we want to show that the process $\hat{Y}$ reaches $0$ with a small probability.

\blemma\label{t4.3}
Suppose that  Condition \ref{c} (iib)-(iic) and condition \eqref{1.4} hold. For the constant $\delta$ in (\ref{7.5})
and small $\epsilon>1$, there is an integer $n_0>0$ so that
 \beqnn
\mbf{P}\{\tau_0^{\hat{Y}}=\infty\}
\geq A_0(\epsilon):=\prod_{n=n_0}^\infty(1-2\epsilon^{\delta(2n-1)}-\epsilon^{2n-3}).
 \eeqnn
\elemma
\proof
In this proof we use $\mbf{E}_{\bar{\varepsilon}}$
and $\mbf{P}_{\bar{\varepsilon}}$ to denote the
conditional
expectation and conditional probability
with respect to $\mathscr{F}_{\hat{Y}(\tau_{\bar{\varepsilon}}^{\hat{Y}})}$.
We first estimate $\mbf{P}_{\epsilon_n}\{\gamma_{n+1}> \epsilon_{n+1}^{-l}|\gamma_n>\epsilon_n^{-l}\}$.
Recall that $\sigma^{Y_{n+1}}(\epsilon_{n-1}):=\inf\{t\ge0: Y_{n+1}(t)>\epsilon_{n-1}\}$.
By Fatou's lemma and \eqref{7.3},
 \beqnn
 \ar\ar
\epsilon_{n-1}
\mbf{P}_{\epsilon_n}\big\{\sigma^{Y_{n+1}}(\epsilon_{n-1})<\gamma_{n+1}
\big|\gamma_n>\epsilon_n^{-l}\big\} \cr
 \ar\ar\qquad\le
\mbf{E}_{\epsilon_n}
\Big[Y_{n+1}(\sigma^{Y_{n+1}}(\epsilon_{n-1})\wedge\gamma_{n+1})
1_{\{\sigma^{Y_{n+1}}(\epsilon_{n-1})<\gamma_{n+1}\}}
\big|\gamma_n>\epsilon_n^{-l}\Big] \cr
 \ar\ar\qquad\le
\mbf{E}_{\epsilon_n}
\Big[Y_{n+1}(\sigma^{Y_{n+1}}(\epsilon_{n-1})\wedge\gamma_{n+1})
\big|\gamma_n>\epsilon_n^{-l}\Big] \cr
 \ar\ar\qquad\le
\liminf_{t\to\infty}
\mbf{E}_{\epsilon_n}
\Big[Y_{n+1}(t\wedge\sigma^{Y_{n+1}}(\epsilon_{n-1})\wedge\gamma_{n+1})
\big|\gamma_n>\epsilon_n^{-l}\Big]
\le\epsilon_n,
 \eeqnn
which implies
 \beqlb\label{5.4}
\mbf{P}_{\epsilon_n}\big\{\sigma^{Y_{n+1}}(\epsilon_{n-1})<\gamma_{n+1}
\big|\gamma_n>\epsilon_n^{-l}\big\}
 \leq
\epsilon_{n}/\epsilon_{n-1}=\epsilon^{2n-1}.
 \eeqlb
Note that by \eqref{7.12},
 \beqlb\label{7.4}
\hat{X}(T_n)\le\gamma_n^{-\frac{1}{p +\delta}}<\epsilon_n^\frac{l}{p +\delta}\quad\mbox{ given  } \gamma_n>\epsilon_n^{-l}
 \eeqlb
for all $n\ge1$.

For $\delta>0$,
by \eqref{7.3}, the definition of the process
$(Y_n(t))_{t\ge0}$ and It\^o's formula,
with respect to  $\{\mathscr{F}_{Y_n(\tau_{\epsilon_n}^{Y_n})}\}$
and for $\gamma_n>\epsilon_n^{-l}$,
 \beqnn
t\mapsto Y_{n+1}(t\wedge\gamma_{n+1})^{-\delta} \exp\Big\{-\int_0^{t\wedge\gamma_{n+1}}
G_\delta(\hat{X}(T_n), Y_{n+1}(s))\dd s\Big\}
 \eeqnn
is a martingale,
where
 \beqnn
G_\delta(u,v)
 :=
\delta[b_1(v)+\theta(v)u^\kappa]v^{-1}
+\frac{\delta(\delta+1)}{2}b_2(v)v^{-2}
+\delta(\delta+1)b_3(v)H_{2,\delta}(v)
 \eeqnn
with the function $H_{2,\delta}$ defined in \eqref{0.4}.
Taking an expectation and using Fatou's lemma, for all $n\ge1$, we have
 \beqlb\label{7.6}
\epsilon_n^{-\delta}
 \ar=\ar
\liminf_{t\to\infty}\mbf{E}_{\epsilon_n}\Big[Y_{n+1}(t\wedge\gamma_{n+1})^{-\delta} \exp\Big\{-\int_0^{t\wedge\gamma_{n+1}}
G_\delta(\hat{X}(T_n),Y_{n+1}(s))\dd s\Big\}\big|\gamma_n>\epsilon_n^{-l}\Big] \cr
 \ar\ge\ar
\mbf{E}_{\epsilon_n}\Big[\liminf_{t\to\infty}Y_{n+1}(t\wedge\gamma_{n+1})^{-\delta} \exp\Big\{-\int_0^{t\wedge\gamma_{n+1}}
G_\delta(\hat{X}(T_n),Y_{n+1}(s))\dd s\Big\}\big|\gamma_n>\epsilon_n^{-l}\Big] \cr
 \ar=\ar
\epsilon_{n+1}^{-\delta}\mbf{E}_{\epsilon_n}\Big[ \exp\Big\{-\int_0^{\gamma_{n+1}} G_\delta(\hat{X}(T_n),Y_{n+1}(s))\dd
s\Big\}\big|\gamma_n>\epsilon_n^{-l}\Big].
 \eeqlb

Observe that under \eqref{7.5},
there is a constant $n_0>1$ so that for all $n\ge n_0$,
 \beqlb\label{7.9}
\delta_0(n)
 \ar:=\ar
(n+1)^2(-l+\theta-1)+n^2\frac{l\kappa}{p +\delta} \cr
 \ar=\ar
(n+1)^2[\frac{l\kappa}{p +\delta}-l+\theta-1]
-(2n+1)\frac{l\kappa}{p +\delta}>0
 \eeqlb
and
 \beqlb\label{7.10}
\delta_1(n):=
-l(n+1)^2+q(n-1)^2
=(q-l)(n-1)^2-4nl>0.
 \eeqlb

Under Condition \ref{c} (iib) and (iic), for all $0<u,v\le c^*$, we have
 \beqnn
G_\delta(u,v)
\le
\delta c_\theta u^\kappa v^{\theta-1}
+
\delta(\delta+1)G_{2,0}(v)
\le
\delta(\delta+1)(c_\theta\vee b)[u^\kappa v^{\theta-1}
+v^q]
 \eeqnn
with $G_{2,0}$ given in \eqref{0.8} and
 \beqnn
\epsilon_{n+1}\le Y_{n+1}(s)\le \epsilon_{n-1}
\quad\text{for}\quad s< \gamma_{n+1}\wedge \sigma^{Y_{n+1}}(\epsilon_{n-1}).
 \eeqnn
Then by \eqref{7.4}, given
$\gamma_n>\epsilon_n^{-l}$ and for $s< \gamma_{n+1}\wedge \sigma^{Y_{n+1}}(\epsilon_{n-1})$,
 \beqnn
G_\delta(\hat{X}(T_n),Y_{n+1}(s))
 \ar\le\ar
\delta(\delta+1)(c_\theta\vee b)
[\hat{X}(T_n)^\kappa Y_{n+1}(s)^{\theta-1}
+Y_{n+1}(s)^q] \cr
 \ar\le\ar
\delta(\delta+1)(c_\theta\vee b)
[\epsilon_n^\frac{\kappa l}{p +\delta}\epsilon_{n+1}^{\theta-1}
+\epsilon_{n-1}^q].
 \eeqnn
It follows from \eqref{7.9} and \eqref{7.10} that, given
$\gamma_n>\epsilon_n^{-l}$ and $\gamma_{n+1}<\epsilon_{n+1}^{-l}\wedge \sigma^{Y_{n+1}}(\epsilon_{n-1})$,
 \beqlb\label{7.11}
 \ar\ar
\int_0^{\gamma_{n+1}}G_\delta(\hat{X}(T_n),Y_{n+1}(s))\dd s \cr
 \ar\ar\quad\le
\delta(\delta+1)(c_\theta\vee b)\gamma_{n+1}
[\epsilon_n^\frac{\kappa l}{p +\delta}\epsilon_{n+1}^{\theta-1}
+\epsilon_{n-1}^q] \cr
 \ar\ar\quad\le
\delta(\delta+1)(c_\theta\vee b) \epsilon_{n+1}^{-l}\big(
\epsilon_{n}^{\frac{\kappa l}{p +\delta}}\epsilon_{n+1}^{\theta-1}
+\epsilon_{n-1}^{q}\big) \cr
 \ar\ar\quad=
\delta(\delta+1)(c_\theta\vee b)
[\epsilon^{\delta_1(n)}+\epsilon^{\delta_0(n)}]\le \ln 2
 \eeqlb
for all $n\ge n_0$ and small enough $\epsilon$.
From \eqref{7.6} and \eqref{7.11}  it follows that
all $n\ge n_0$ and small enough $\epsilon$,
 \beqnn
\epsilon_n^{-\delta}
 \ar\geq\ar
\epsilon_{n+1}^{-\delta}\mbf{E}_{\epsilon_n}
\Big[\exp\Big\{-\int_0^{\gamma_{n+1}} G_\delta(\hat{X}(T_n),Y_{n+1}(s))\dd s\Big\} 1_{\{\gamma_{n+1}<\epsilon_{n+1}^{-l}\wedge
\sigma^{Y_{n+1}}(\epsilon_{n-1})\}}
\big|\gamma_n>\epsilon_n^{-l}\Big]  \cr
 \ar\geq\ar
2^{-1}\epsilon_{n+1}^{-\delta}  \mbf{P}_{\epsilon_n}\big\{\gamma_{n+1}<\epsilon_{n+1}^{-l}\wedge
\sigma^{Y_{n+1}}(\epsilon_{n-1})\big|\gamma_n>\epsilon_n^{-l}\big \} \cr
 \ar\geq\ar
2^{-1}\epsilon_{n+1}^{-\delta} \Big[\mbf{P}_{\epsilon_n}\big\{\gamma_{n+1}<\epsilon_{n+1}^{-l}
\big|\gamma_n>\epsilon_n^{-l}\big\}
-\mbf{P}_{\epsilon_n}\big\{\gamma_{n+1}> \sigma^{Y_{n+1}}(\epsilon_{n-1}) \big|\gamma_n>\epsilon_n^{-l}\big\}\Big].
 \eeqnn
It follows from \eqref{5.4} that
for all $n\ge n_0$ and small enough $\epsilon$,
 \beqnn
\mbf{P}_{\epsilon_n}\big\{\gamma_{n+1}<\epsilon_{n+1}^{-l}
\big|\gamma_n>\epsilon_n^{-l}\big\}
 \ar\le\ar
2\epsilon_n^{-\delta}\epsilon_{n+1}^{\delta}
+\mbf{P}_{\epsilon^n}\big\{\gamma_{n+1}>\sigma^{Y_{n+1}}(\epsilon_{n-1})
\big|\gamma_n>\epsilon_n^{-l}\big\} \cr
 \ar\le\ar
{2}\epsilon^{\delta(2n+1)}+\epsilon^{2n-1}.
 \eeqnn

Observe that by the Markov property,
 \beqnn
\mbf{P}\{\cap_{n=n_0}^m\{\gamma_{n}>\epsilon_{n}^{-l}\} \}
 \ar=\ar
\mbf{E}\Big[\mbf{P}_{\epsilon_{m-1}}\big\{\gamma_{m}>\epsilon_{m}^{-l}| \cap_{n=n_0}^{m-1}\{\gamma_{n}>\epsilon_{n}^{-l}\}\big\}
1_{\cap_{n=n_0}^{m-1}\{\gamma_{n}>\epsilon_{n}^{-l}\}}\Big] \cr
 \ar=\ar
\mbf{E}\Big[\mbf{P}_{\epsilon_{m-1}}\big\{\gamma_{m}>\epsilon_{m}^{-l}| \gamma_{m-1}>\epsilon_{m-1}^{-l}\big\}
1_{\cap_{n=n_0}^{m-1}\{\gamma_{n}>\epsilon_{n}^{-l}\}}\Big] \cr
 \ar\geq\ar
(1-2\epsilon^{\delta(2m-1)}-\epsilon^{2m-3}) \mbf{P}\{\cap_{n=n_0}^{m-1}\{\gamma_{n}>\epsilon_{n}^{-l}\} \} \cr
 \ar\geq\ar
\prod_{n=n_0}^m(1-2\epsilon^{\delta(2n-1)}-\epsilon^{2n-3}).
 \eeqnn
Letting $m\to\infty$ we get
 \beqnn
A_0(\epsilon)\le\mbf{P}\big\{\cap_{n=n_0}^\infty\{\gamma_{n}>\epsilon_{n}^{-l}\} \big\}
\le
\mbf{P}\{\tau_0^{\hat{Y}}=\infty\},
 \eeqnn
which ends the proof.
\qed

\blemma\label{t4.5}
Under Condition \ref{c} (iia) with $p>0$ and (iic) and (iii), for each $\delta>0$  and small enough $\varepsilon>0$,
there are constants $C(\delta,\varepsilon)>0$ and $\delta_1\in(0,1)$ that do not
depend on $\varepsilon$ so that for all $X_0=\varepsilon^m$
with large enough $m$ we have
 \beqlb\label{5.9}
\mbf{P}\{Y_t\ge \hat{Y}(t)~\mbox{for all }~ t\ge0\}=:\mbf{P}\{B\}
\ge1-C(\delta,\varepsilon)\varepsilon^{m\delta_1/8}.
 \eeqlb
\elemma
\proof
By Lemma \ref{t4.2}, there are constants $C(\delta,\varepsilon)>0$
and $\delta_1\in(0,1)$ so that for all $X_0=\varepsilon^m$
with large $m$ we have
 \beqnn
\mbf{P}\{X_t\le \hat{X}(t)~\mbox{for all}~ t\ge0\}=:\mbf{P}(A)\ge1-C(\delta,\varepsilon)\varepsilon^{m\delta_1/8}.
 \eeqnn

Observe that under Condition \ref{c}(iic),
given $A$ and $s\ge T$, we have
 \beqnn
\kappa(X_s)
\le X_s^\kappa
\le \hat{X}(s)^\kappa
\le \hat{X}(T)^\kappa.
 \eeqnn
Since $(B_t)_{t\ge0}$,
$(W_t)_{t\ge0}$,
$\{\tilde{M}(\dd t,\dd z , \dd u)\}$ and
$\{\tilde{N}(\dd t,\dd z , \dd u)\}$ are independent,
then by using \eqref{7.3}, the definition of $(\hat{Y}(t))_{t\ge0}$,
Condition \ref{c}(iii) and Proposition \ref{t4.4}, $\mbf{P}\{B|A\}=1$.
It follows that
 \beqnn
\mbf{P}\{B\}
\ge\mbf{P}\{B|A\}\cdot
\mbf{P}\{A\}
\ge1-C(\delta,\varepsilon)\varepsilon^{m\delta_1/8},
 \eeqnn
which ends the proof.
\qed

\blemma\label{t4.6}
Under Condition \ref{c} (iii), for each $\varepsilon>0$, there is a constant $t_0>0$ so that
 \beqlb\label{4.11}
\mbf{P}\Big\{\sup_{t\ge t_0}X_t\le \varepsilon,Y_{t_0}>0\Big\}>0.
 \eeqlb
\elemma
\proof
Since $X_t\rightarrow 0$ as $t\rightarrow \infty$ by Lemma \ref{t2.1} (i),
there are constants $t_0>0$ and $n\ge1$ so that
 \beqlb\label{3.16}
\mbf{P}\Big\{\sup_{t\ge t_0}X_t\le \varepsilon,
\sup_{0\le t<t_0}X_t\le n \Big\}>0.
 \eeqlb
Let $(\tilde{Y}_t)_{t\ge 0}$ be the nonnegative solution to
 \begin{equation}\label{3.18}
 	\begin{split}
\tilde{Y}_t
=&
Y_0-\int_0^t[b_1(\tilde{Y}_s)+C_n\theta(\tilde{Y}_s)]\dd s
+\int_0^t b_2(\tilde{Y}_s)^{1/2}\dd W_s \\
&+\int_0^t\int_0^\infty\int_0^{b_3(\tilde{Y}_{s-})}z \tilde{N}(\dd s,\dd z,\dd u),
 \end{split}
\end{equation}
where $C_n:=\sup_{x\in[0,n]}\kappa(x)$.
 Under Condition \ref{c} (iii), by the comparison theorem (Proposition \ref{t4.4}),
 \beqlb\label{3.17}
\mbf{P}\Big\{Y_t\ge \tilde{Y}_t
\mbox{  for all } 0\le t\le t_0\big|\sup_{0\le t\le t_0}X_t\le n\Big\}=1.
 \eeqlb
It is easy to see that
$\mbf{P}\{\tilde{Y}_{t_0}>0\}>0$.
Since $(\tilde{Y}_t)_{t\ge0}$ and $(X_t)_{t\ge0}$ are independent,
then by \eqref{3.16} we get
 \beqnn
\mbf{P}\Big\{\sup_{t\ge t_0}X_t\le \varepsilon,
\sup_{0\le t<t_0}X_t\le n,\tilde{Y}_{t_0}>0\Big\}
=\mbf{P}\Big\{\sup_{t\ge t_0}X_t\le \varepsilon,
\sup_{0\le t<t_0}X_t\le n\Big\}
\mbf{P}\{\tilde{Y}_{t_0}>0\}
>0.
 \eeqnn
Using \eqref{3.17} we get
 \beqnn
\mbf{P}\Big\{\sup_{t\ge t_0}X_t\le \varepsilon,
\sup_{0\le t<t_0}X_t\le n,Y_{t_0}>0\Big\}>0,
 \eeqnn
which implies \eqref{4.11}.
\qed

\noindent{\it Proof of Theorem \ref{t0.02}.}
We first show that
for given $X_0=\varepsilon^m$ and $Y_0$ with $m$ large enough
and $\varepsilon$ small enough,
there is a constant $C(\varepsilon)>0$
so that
 \beqlb\label{7.13}
\mbf{P}\{\tau_0^Y=\infty\}\ge C(\varepsilon).
 \eeqlb

Let $B^c$ denote the complementary set of $B$, which is given in \eqref{5.9}.
By Lemma \ref{t4.5},
there are constants $C(\delta,\varepsilon)>0$ and $\delta_1\in(0,1)$
%independent of $\delta$
independent of $\delta$ so that
 \beqlb\label{7.8}
\mbf{P}\{B^c\}
\le C(\delta,\varepsilon)\varepsilon^{m\delta_1/8}.
 \eeqlb
Observe that
 \beqnn
\{\tau_0^{\hat{Y}}=\infty\}
 \ar=\ar
(\{\tau_0^{\hat{Y}}=\infty\}\cap B)
\cup
(\{\tau_0^{\hat{Y}}=\infty\}\cap B^c) \cr
 \ar\subset\ar
(\{\tau_0^Y=\infty\}\cap B)
\cup B^c \subset
\{\tau_0^Y=\infty\}
\cup
B^c.
 \eeqnn
Therefore, by Lemma \ref{t4.3} and \eqref{7.8},
 \beqnn
\mbf{P}\{\tau_0^Y=\infty \}
 \ar\ge\ar
\mbf{P}\{\{\tau_0^Y=\infty\}\cup B^c\}
-\mbf{P}\{B^c\} \cr
 \ar\ge\ar
\mbf{P}\{\tau_0^{\hat{Y}}=\infty\}
-\mbf{P}\{B^c\}
\ge A_0(\epsilon)
-C(\delta,\varepsilon)\varepsilon^{m\delta_1/8}>0
 \eeqnn
for $m$ large enough and small enough $\epsilon$
and $\varepsilon$, which gives \eqref{7.13} for some constant
$C(\varepsilon)>0$.

By Lemma \ref{t4.6}, for each $\varepsilon>0$, there is a constant $t_0:=t_0(\varepsilon)>0$ so that
$\mbf{P}\{X_{t_0}\le \varepsilon,Y_{t_0}>0\}>0$.
By the Markov property and \eqref{7.13}, for each $t>0$ and small
enough $\varepsilon>0$, there is a constant $C(\varepsilon)>0$
so that for $X_t\le \varepsilon$ and $Y_t>0$, we have
 \beqnn
\mbf{P}\{\tau_0^Y=\infty|(X_t,Y_t)\}
\ge C(\varepsilon).
 \eeqnn
 It follows that
 \beqnn
\mbf{P}\big\{\tau_0^Y=\infty\big\}
=\mbf{P}\big\{X_{t_0}\le \varepsilon,Y_{t_0}>0\big\}
\cdot
\mbf{P}\big\{\tau_0^Y=\infty
|X_{t_0}\le \varepsilon,Y_{t_0}>0\big\}>0,
 \eeqnn
which ends the proof.  \qed

\subsection{Proof of Theorem \ref{t0.06}}\label{3.5}

For $\delta\in(-1,0)\cup(0,\infty)$, and $i=1,2$ recall the definitions of $H_{i,\delta},H_{i,0},G_{i,0}$ in \eqref{0.3}--\eqref{0.8}.
For $x,y,\beta>0$ and $r\in(-(\beta^{-1}\wedge1),1)$ define
 \beqlb\label{4.12}
G_r(x,y):=\beta G_{1,r}(x)-G_{2,r}(y)-\kappa(x)\theta(y)y^{-1}
 \eeqlb
with
 \beqlb\label{4.12a}
G_{1,r}(x):=a_1(x)x^{-1}+ (1+\beta r)2^{-1}a_2(x)x^{-2}+(1+\beta r)a_3(x)H_{1,\beta r}(x)
 \eeqlb
and
 \beqlb\label{4.12b}
G_{2,r}(y):=b_1(y)y^{-1}
+(1-r)2^{-1}b_2(y)y^{-2}+(1-r)b_3(y)H_{2,-r}(y).
 \eeqlb
Thus,
 \beqlb\label{4.12c}
G_0(x,y):=\beta G_{1,0}(x)-G_{2,0}(y)-\kappa(x)\theta(y)y^{-1}.
 \eeqlb
To prove Theorem \ref{t0.06}, we first prove  the following assertions.

\blemma\label{t4.7}
For any $x,y>0$, we have
 \beqlb\label{3.15a}
G_r(x,y)
\le
(1-r)G_0(x,y)+\beta(\beta+1)rG_{1,0}(x),\qquad r\in (0,1)
 \eeqlb
and
 \beqlb\label{3.15b}
G_r(x,y)
\ge
(1+\beta r)G_0(x,y)+(\beta+1)rG_{2,0}(y)
+\beta r\kappa(x)\theta(y)y^{-1},~ r\in(-(\beta^{-1}\wedge1),0).
 \eeqlb
\elemma
\proof
Observe that for each $i=1,2$, $x>0$ and $r\ge0$, we have
 \beqnn
H_{i,r}(x)\le H_{i,0}(x),\quad H_{i,-r}(x)\ge H_{i,0}(x),
 \eeqnn
which implies that
 \beqnn
G_{1,r}(x)\le (1+\beta r) G_{1,0}(x),\quad
G_{2,r}(x)\ge (1-r) G_{2,0}(x).
 \eeqnn
Then for $x,y>0$,
 \beqnn
G_r(x,y)
 \ar\le\ar
\beta (1+\beta r)G_{1,0}(x)-(1-r)G_{2,0}(y)
-(1-r)\kappa(x)\theta(y)y^{-1} \cr
 \ar=\ar
(1-r)G_0(x,y)+\beta(\beta+1)rG_{1,0}(x),~~ r\in (0,1)
 \eeqnn
and
 \beqnn
G_r(x,y)
 \ar\ge\ar
\beta (1+\beta r)G_{1,0}(x)-(1-r)G_{2,0}(y)
-\kappa(x)\theta(y)y^{-1} \cr
 \ar=\ar
(1+\beta r)G_0(x,y)+(\beta+1)rG_{2,0}(y)
+\beta r\kappa(x)\theta(y)y^{-1},~~r\in(-(\beta^{-1}\wedge1),0),
 \eeqnn
which finishes the proof.
\qed

The following result is key to the proof of Theorem \ref{t0.06}.

\blemma\label{t3.1}
Under the assumptions of Theorem \ref{t0.06},
for any $0<\varepsilon_1<\varepsilon$, if $X_0,Y_0\le\varepsilon_1$, then
 we have
 \beqlb\label{4.10}
\mbf{P}\{\tau_0^Y\wedge\sigma_\varepsilon^{X}
\wedge\sigma_\varepsilon^{Y}<\infty\}=1.
 \eeqlb
\elemma
\proof
The proof is an application of Corollary \ref{cor2.1}.
We first present the key function $g$ satisfying the conditions of Corollary \ref{cor2.1}.
Define $g(u):=\e^{-\lambda u^r}$ for $u,\lambda>0$ and $0<r<1$.
Let $0<\varepsilon<c^*$ (determined in Condition \ref{c}(i)) and $g(x,y):=g(x^{-\beta}y)$ for all $x,y>0$,
where the value of constant $\beta>0$ is to be specified later.
In the following we show that there are constants $d_1,d_2>0$ so that for all
$0<x,y<\varepsilon$, we have, respectively,
\begin{equation}\label{6.13}
\begin{split}
&
Lg(x,y)\ge r\lambda d_1 g(x,y)
\mbox{ under condition (i) of Theorem \ref{t0.06}} \\
&
\mbox{and }
Lg(x,y)\ge r\lambda d_2x^p g(x,y)
\mbox{  under condition (ii) of Theorem \ref{t0.06}}.
\end{split}
\end{equation}

Recall the definitions of $K_z^1$ and $K_z^2$ in (\ref{3.1})
and $G_{r},G_{1,r},G_{2,r}$ in \eqref{4.12}-\eqref{4.12b}.
For simplicity we denote  $u=x^{-\beta}y$ in the following.
By \eqref{3.1} and \eqref{7.7}-\eqref{7.7b},
 \beqnn
g(u)^{-1}K_z^1g(x,y)
 \ge
-r\beta(r\beta+1)\lambda u^rz^2x^{-2}\int_0^1(1+z x^{-1}v)^{-r\beta-2}(1-v)\dd v
 \eeqnn
and
 \beqnn
g(u)^{-1} K_z^2g(x,y)
 \ge
r(1-r)\lambda u^rz^2y^{-2}\int_0^1(1+z y^{-1}v)^{r-2}(1-v)\dd v.
 \eeqnn
Then one can get
 \beqnn
L_1g(x,y)
 \ar=\ar
-\lambda r\beta g(u)u^r a_1(x)/x
+2^{-1}\big[(\lambda r\beta)^2g(u)u^{2r} \cr
 \ar\ar
-\lambda r\beta(1+r\beta)g(u)u^r\big]a_2(x)/x^2
+a_3(x)\int_0^\infty
K_z^1g(x,y) \mu(\dd z) \cr
 \ar\ge\ar
-\lambda r\beta g(u)u^rG_{1,r}(x)
 \eeqnn
and
 \beqnn
L_2g(x,y)
 \ar=\ar
\lambda r g(u)u^r [b_1(y)+\kappa(x)\theta(y)]/y
+2^{-1}\big[(\lambda r)^2g(u)u^{2r} \cr
 \ar\ar
+\lambda r(1-r)g(u)u^r\big]b_2(y)/y^2
+b_3(x)\int_0^\infty
K_z^2g(x,y) \nu(\dd z) \cr
 \ar\ge\ar
\lambda r g(u)u^r[G_{2,r}(y)+\kappa(x)\theta(y)y^{-1}].
 \eeqnn
Thus,
 \beqlb\label{6.8}
Lg(x,y)
=
L_1g(x,y)+L_2g(x,y)
\ge
-\lambda r  u^rg(u)G_r(x,y)=-\lambda r u^rG_r(x,y)g(x,y).
 \eeqlb

In the following we use the inequality \eqref{3.15a} in Lemma \ref{t4.7}
to estimate $G_r$. Recall the definition of $G_0$ in \eqref{4.12c}.

Under condition (i) of Theorem \ref{t0.06}, taking $\beta=\kappa/(1-\theta)$,
we have $-\beta a+b>0$.
Then there exist small constants $c_1>0$ and $0<r<1-\theta$ so that
 \beqlb\label{7.21}
(1-r)(-\beta a+b)-a\beta(\beta+1)r \ge c_1
 \eeqlb
Under condition (i) of Theorem \ref{t0.06} we have $p=q=0$ and then using Condition \ref{c} (i) we obtain
 \beqlb\label{7.22}
-G_0(x,y)
 \ar=\ar
-\beta G_{1,0}(x)+G_{2,0}(y)+\kappa(x)\theta(y)y^{-1} \cr
 \ar\ge\ar
-\beta a+b+c_\theta x^\kappa y^{\theta-1}
=
-\beta a+b+c_\theta u^{\theta-1}
 \eeqlb
and $G_{1,0}(x)\le a$ for all $0<x,y<\varepsilon$.
It thus follows from \eqref{3.15a} and \eqref{7.21}--\eqref{7.22} that
 \beqnn
-G_r(x,y)
 \ar\ge\ar
-(1-r)G_0(x,y)-\beta(\beta+1)rG_{1,0}(x) \cr
 \ar\ge\ar
(1-r)(-\beta a+b)+(1-r)c_\theta u^{\theta-1}-a\beta(\beta+1)r  \cr
 \ar\ge\ar
c_1+(1-r)c_\theta u^{\theta-1}
 \eeqnn
for all $0<x,y<\varepsilon$. It then follows from Lemma \ref{t2.3} that
 \beqnn
-u^rG_r(x,y)\ge c_1u^r+(1-r)c_\theta u^{r+\theta-1}\ge d_1,\qquad 0<x,y<\varepsilon
 \eeqnn
for some constant $d_1>0$.
Then the first part of \eqref{6.13} follows from \eqref{6.8}.

Under condition (ii) of Theorem \ref{t0.06} and Condition \ref{c} (i),
taking $\beta=p/q$,
we have $p=\beta q=\kappa-\beta(1-\theta)$ and then
 \beqlb\label{4.9}
-G_0(x,y)
 \ar\ge\ar
-\beta ax^p+by^q +c_\theta x^\kappa y^{\theta-1} \cr
 \ar=\ar
x^p[-\beta a+bx^{-p}y^q+c_\theta x^{\kappa-p} y^{\theta-1}] \cr
 \ar=\ar
x^p[-\beta a+bu^q+c_\theta u^{\theta-1}],\quad 0<x,y<\varepsilon.
 \eeqlb
For
 \beqnn
\bar{p}:=1+q/(1-\theta)=[q+(1-\theta)](1-\theta)^{-1},
\quad\bar{q}:
=\bar{p}/(\bar{p}-1)
=[q+(1-\theta)]q^{-1},
 \eeqnn
we have $q/\bar{p}+(\theta-1)/\bar{q}=0$
and then by Lemma \ref{t2.3},
 \beqlb\label{4.4}
bu^q+c_\theta u^{\theta-1}
 \ar\ge\ar
\bar{p}^{1/\bar{p}}
\bar{q}^{1/\bar{q}}
b^{1/\bar{p}}c_\theta^{1/\bar{q}}
u^{q/\bar{p}+(\theta-1)/\bar{q}}
=
\bar{p}^{1/\bar{p}}
\bar{q}^{1/\bar{q}}
b^{1/\bar{p}}c_\theta^{1/\bar{q}} \cr
 \ar=\ar
[q+(1-\theta)]\Big(\frac{b}{1-\theta}\Big)^{\frac{1-\theta}{q+1-\theta}}
\Big(\frac{c_\theta}{q}\Big)^{\frac{q}{q+1-\theta}}=:c_2.
 \eeqlb
Under condition \eqref{1.6}, we have $c_2>\beta a$.
It follows from \eqref{4.9} that
 \beqnn
-G_0(x,y)x^{-p}\ge-\beta a+bu^q+c_\theta u^{\theta-1}
\ge c_2-\beta a>0,\quad 0<x,y<\varepsilon.
 \eeqnn
Then by \eqref{3.15a} and Condition \ref{c} (i) again, there are constants $0<r<1-\theta$ and
$c_3:=(c_2-\beta a)(1-r)$ so that
$c_3>r\beta(\beta+1)a$
and
 \beqlb\label{4.3}
-G_r(x,y)x^{-p}
\ge
-(1-r)G_0(x,y)x^{-p}-r\beta(\beta+1)G_{1,0}(x)x^{-p}
\ge c_3-r\beta(\beta+1)a>0
 \eeqlb
for all $0<x,y<\varepsilon$.
Then
 \beqlb\label{7.15}
-G_r(x,y)u^r
\ge [c_3-r\beta(\beta+1)a]x^p>0,
\qquad 0<x,y<\varepsilon,\,u\ge 1.
 \eeqlb
Let $\delta>0$ be a constant satisfying
 \beqlb\label{7.14}
((1-\delta)c_2-\beta a)(1-r)
-r\beta(\beta+1)a>0.
 \eeqlb
By \eqref{4.9} and \eqref{4.4} we obtain
 \beqnn
-G_0(x,y)
 \ar\ge\ar
x^p\big[-\beta a+(1-\delta)(bu^q+c_\theta u^{\theta-1})
+\delta(bu^q+c_\theta u^{\theta-1})\big] \cr
 \ar\ge\ar
x^p\big[((1-\delta)c_2-\beta a)+\delta c_\theta u^{\theta-1}\big],
~~ 0<x,y<\varepsilon
 \eeqnn
and then by \eqref{7.14} and the same argument as in \eqref{4.3},
 \beqnn
-G_r(x,y)u^r
 \ar\ge\ar
x^p u^r\big[((1-\delta)c_2-\beta a)(1-r)
-r\beta(\beta+1)a
+\delta c_\theta u^{\theta-1}\big] \cr
 \ar\ge\ar
\delta c_\theta x^p u^{r+\theta-1}
\ge \delta c_\theta x^p,
\qquad 0<x,y<\varepsilon,\,u\le 1.
 \eeqnn
This and \eqref{7.15} imply the second part of \eqref{6.13} by \eqref{6.8}.

%{\bf Step 2.}
Letting $\kappa\bar{p}=p$ in Lemma \ref{t1.5}, we have
$\int_0^\infty X_s^p\dd s=\infty$
almost surely.
Since $X_t\to0$ as $t\to\infty$ by Lemma \ref{t2.1} (i),
then for all $0<\varepsilon<c^*$, we have $\int_0^\infty (X_s^p\wedge \varepsilon)\dd s=\infty$
almost surely.
Using the above assertions and
Corollary \ref{cor2.1},
we have
 \beqnn
\mbf{P}\{\tau_0^X\wedge\tau_0^Y
\wedge\sigma_\varepsilon^X\wedge\sigma_\varepsilon^Y<\infty
\}
\ge \e^{-\lambda (X_0^{-\beta}Y_0)^r}
 \eeqnn
for all $0<\varepsilon_1<\varepsilon$ and $\lambda>0$.
Since $\tau_0^Y=\infty$ almost surely by Lemma \ref{t2.1} (ii),
then letting $\lambda\to0$ in the above inequality
we have \eqref{4.10}.
\qed

Now we are ready to prove Theorem \ref{t0.06}.

\noindent{\it Proof of Theorem \ref{t0.06}}.
For any $0<\varepsilon_1<\varepsilon$, if $X_0,Y_0\le\varepsilon_1$,
by Lemma \ref{t3.2} we have
\beqnn
\mbf{P}\{\sigma_{\varepsilon}^X<\infty\}
=\mbf{P}\big\{\sup_{t\ge0}X_t\ge\varepsilon\big\}
\le C(\varepsilon_1/\varepsilon)^{1/4}\quad\text{and}\quad
\mbf{P}\{\sigma_{\varepsilon}^Y<\infty\}
\le C(\varepsilon_1/\varepsilon)^{1/4}
\eeqnn
for some constant $C>0$,
which implies
\beqnn
\mbf{P}\{\sigma_\varepsilon^{X}\wedge\sigma_\varepsilon^{Y}<\infty\}
\le
\mbf{P}\{\sigma_{\varepsilon}^X<\infty\}
+\mbf{P}\{\sigma_{\varepsilon}^Y<\infty\}
\le 2C(\varepsilon_1/\varepsilon)^{1/4}.
\eeqnn
It follows from Lemma \ref{t3.1} that
\beqlb\label{proof1.6}
\mbf{P}\{\tau_0^Y<\infty\}
\ge
\mbf{P}\{\tau_0^Y\wedge\sigma_\varepsilon^{X}\wedge\sigma_\varepsilon^{Y}<\infty\}
-\mbf{P}\{\sigma_\varepsilon^{X}\wedge\sigma_\varepsilon^{Y}<\infty\}
\ge
1-2C(\varepsilon_1/\varepsilon)^{1/4}.
\eeqlb
If $X_0>\varepsilon_1 $ or $Y_0>\varepsilon_1$, (\ref{proof1.6}) also holds
following from  Lemma \ref{t2.1} (i) and the Markov property.
Letting $\varepsilon_1\to0$ we finish the proof.

\bremark
Under Condition \ref{c} (i) and condition \eqref{1.5},
for $q>0$, take $\beta=p/q$
and $\varepsilon>0$ small enough so that
 \beqnn
\frac{ap}{q(q+1-\theta)}
<\Big(\frac{b}{1-\theta}\Big)^{\frac{1-\theta}{q+1-\theta}}
\cdot
\Big(\frac{c_\theta \varepsilon^{-\theta_1}}{q}\Big)^{\frac{q}{q+1-\theta}},
 \eeqnn
where $\theta_1:=\beta(1-\theta)-(\kappa-p)>0$.
Then by an argument similar to that for \eqref{4.9}, we get
 \beqnn
-G_0(x,y)
 \ar\ge\ar
x^p[-\beta a+bu^q+c_\theta \varepsilon^{-\theta_1} u^{\theta-1}],\quad 0<x,y<\varepsilon.
 \eeqnn
By essentially the same argument after \eqref{4.9} in the proof of
Lemma \ref{t3.1}, we can also obtain \eqref{4.10} for $q>0$ and condition \eqref{1.5}.
Therefore, the method for the proof of Theorem \ref{t0.06}
also works for Theorem \ref{t0.01} except the case $q=0$.
\eremark

\subsection{Proof of Theorem \ref{t0.05}}

Recall the function $G_\delta$ in \eqref{4.12}.
We first prove Theorem \ref{t0.05}  for small related
initial values of $X_0$ and $Y_0$,
where the key idea, inspired by the proof of Lemma \ref{t3.1}, is to consider ratio $Y_t/X_t^{\beta} $ process with the value of $\beta>0$ properly selected using the conditions of Theorem \ref{t0.05}.
We then formulate an exponential martingale
that is similar to those  in \cite{LYZh}, and use the martingale to obtain the desired estimate.

\blemma\label{t3.3}
Under the conditions of Theorem \ref{t0.05},
 there exist constants $\beta>0$
and small $\varepsilon>0$ so that for $X_0\le \varepsilon,Y_0= \varepsilon^\beta$, we have
$\mbf{P}\{\tau_0^Y=\infty\}>0$.
\elemma
\proof
First note that for any $\varepsilon_2>0$ and $\varepsilon=\varepsilon_2^{1+\beta^{-1}}$ for any $\beta>0$,
if $X_0\le \varepsilon,Y_0= \varepsilon^\beta$, then
by Lemma \ref{t3.2} we have
 \beqlb\label{3.3}
\mbf{P}
\Big\{\sup_{s\ge 0}X_s\ge\varepsilon_2\Big\}
+\mbf{P}
\Big\{\sup_{s\ge 0}Y_s\ge\varepsilon_2\Big\}
\le
C[(\varepsilon^\beta\varepsilon_2^{-1})^{1/4}
+(\varepsilon\varepsilon_2^{-1})^{1/4}]
\le
C[\varepsilon_2^{\beta/4}
+\varepsilon_2^{1/(4\beta)}],
 \eeqlb
where $C>0$ is a constant independent of $\varepsilon_2$.

We first show that under conditions
(i) or (ii) of Theorem \ref{t0.05}, there are constants $0<\omega<1$, $\beta>0$ and small enough
$\delta\in(0,1),\,0<\varepsilon_2<c^*$ so that
 \beqlb\label{3.2}
G_{-\delta}(x,y)\ge0,\qquad 0<x,y\le \varepsilon_2,~ x^{-\beta}y\ge \omega
 \eeqlb
and
 \beqlb\label{3.4}
1-\omega^\delta-C[\varepsilon_2^{\beta/4}
+\varepsilon_2^{1/(4\beta)}]>0,
 \eeqlb
where $C>0$ is the constant determined by \eqref{3.3}
and $G_{-\delta}(x,y)$ is defined in \eqref{4.12}.

Under condition (i) of Theorem \ref{t0.05}
and Condition \ref{c} (ii),
select $\beta$ satisfying $b/a<\beta<\kappa/(1-\theta)$.
There are constants $\omega\in(0,1)$ and small enough $\delta>0,\,0<\varepsilon_2<c^*$
so that \eqref{3.4} holds and
 \beqlb\label{7.16}
(1-\beta \delta)[\beta a-b-c_\theta \omega^{\theta-1} \varepsilon_2^{\kappa-\beta(1-\theta)}]-(\beta+1)\delta b
-\beta \delta c_\theta \omega^{\theta-1} \varepsilon_2^{\kappa-\beta(1-\theta)}
>0.
 \eeqlb
Recall function $G_0$ defined in \eqref{4.12c}.
Under condition (i) of Theorem \ref{t0.05}, $p=q=0$
and then by Condition \ref{c} (ii),
for all $0<x,y\le \varepsilon_2$ and $x^{-\beta}y\ge \omega$ we have
 \beqlb\label{7.17}
G_0(x,y)\ge
\beta a-b-c_\theta (x^{-\beta}y)^{\theta-1} x^{\kappa-\beta(1-\theta)}
\ge\beta a-b-c_\theta \omega^{\theta-1} \varepsilon_2^{\kappa-\beta(1-\theta)},
~G_{2,0}(y)\le b
 \eeqlb
and
 \beqlb\label{7.18}
\kappa(x)\theta(y)y^{-1}
\le
c_\theta x^\kappa y^{\theta-1}
=
c_\theta(x^{-\beta}y)^{\theta-1} x^{\kappa-\beta(1-\theta)}
\le
c_\theta \omega^{\theta-1} \varepsilon_2^{\kappa-\beta(1-\theta)}.
 \eeqlb
Moreover, by \eqref{3.15b} in Lemma \ref{t4.7},
\eqref{7.16}-\eqref{7.18}, for all $0<x,y\le \varepsilon_2$ and $x^{-\beta}y\ge \omega$, we have
 \beqnn
G_{-\delta}(x,y)
 \ar\ge\ar
(1-\beta \delta)G_0(x,y)-(\beta+1)\delta G_{2,0}(y)
-\beta \delta\kappa(x)\theta(y)y^{-1} \cr
 \ar\ge\ar
(1-\beta \delta)[\beta a-b-c_\theta \omega^{\theta-1} \varepsilon_2^{\kappa-\beta(1-\theta)}]-(\beta+1)\delta b
-\beta \delta c_\theta \omega^{\theta-1} \varepsilon_2^{\kappa-\beta(1-\theta)}>0,
 \eeqnn
which gives \eqref{3.2}.

Under condition (ii) of Theorem \ref{t0.05} and Condition \ref{c} (ii),
selecting $0<\beta<\kappa/(1-\theta)$, $\omega\in(0,1)$
and small enough $\delta>0,\,0<\varepsilon_2<c^*$ so that \eqref{3.4} holds and
 \beqlb\label{7.19}
(1-\beta \delta)[\beta a-b\varepsilon_2^q-c_\theta \omega^{\theta-1} \varepsilon_2^{\kappa-\beta(1-\theta)}]-(\beta+1)\delta b\varepsilon_2^q
-\beta \delta c_\theta \omega^{\theta-1} \varepsilon_2^{\kappa-\beta(1-\theta)}>0.
 \eeqlb
Since $p=0$, then under Condition \ref{c} (ii), similar to \eqref{7.17},
for all $0<x,y\le  \varepsilon_2$ and $x^{-\beta}y\ge \omega$ we have
 \beqlb\label{7.20}
G_0(x,y)
\ge \beta a-b\varepsilon_2^q-c_\theta \omega^{\theta-1}\varepsilon_2^{\kappa-\beta(1-\theta)},\quad
G_{2,0}(y)\le b \varepsilon_2^q.
 \eeqlb
Thus by \eqref{3.15b} in Lemma \ref{t4.7},
\eqref{7.18}-\eqref{7.20}  and Condition \ref{c} (ii) again, for all $0<x,y\le \varepsilon_2$ and $x^{-\beta}y\ge \omega$, we have
 \beqnn
G_{-\delta}(x,y)
 \ar\ge\ar
(1-\beta \delta)G_0(x,y)-(\beta+1)\delta G_{2,0}(y)
-\beta \delta\kappa(x)\theta(y)y^{-1} \cr
 \ar\ge\ar
(1-\beta \delta)[\beta a-b\varepsilon_2^q-c_\theta \omega^{\theta-1} \varepsilon_2^{\kappa-\beta(1-\theta)}]-(\beta+1)\delta b\varepsilon_2^q
-\beta \delta c_\theta \omega^{\theta-1} \varepsilon_2^{\kappa-\beta(1-\theta)}>0,
 \eeqnn
which proves \eqref{3.2}.

For $v>1$ and $\omega>0$, define stopping times
\[\tau_w:=\inf\{t\ge 0:X_t^{-\beta}Y_t<w\}\quad\text{ and}\quad \sigma_v:=\inf\{t\ge 0:X_t^{-\beta}Y_t>v\},\]
 respectively.
In the following we show that $\mbf{P}\{\tau_w=\infty\}>0$,
which implies the assertion of the lemma.

Let $\omega,\beta,\delta,\varepsilon>0$ be the constants determined in
\eqref{3.2}-\eqref{3.4}.
Let $T:=\tau_w\wedge\sigma_v$.
By \eqref{1.1} and It\^o's formula, for each $\delta>0$,
 \beqnn
(X_{t\wedge T}^{-\beta}Y_{t\wedge T})^{-\delta}\exp\Big\{\delta\int_0^{t\wedge T}G_{-\delta}(X_s,Y_s)\dd s\Big\}
 \eeqnn
is a martingale.
It follows from Fatou's lemma that
 \beqnn
1
 \ar\ge\ar
\mbf{E}[(X_0^{-\beta}Y_0)^{-\delta}]
 =
\liminf_{t\to\infty}\mbf{E}\Big[(X_{t\wedge T}^{-\beta}Y_{t\wedge T})^{-\delta}
\exp\Big\{\delta\int_0^{t\wedge T}G_{-\delta}(X_s,Y_s)\dd s\Big\}\Big]
\cr
 \ar\ge\ar
\mbf{E}\Big[
\liminf_{t\to\infty}(X_{t\wedge T}^{-\beta}Y_{t\wedge T})^{-\delta}\exp\Big\{\delta\int_0^{t\wedge T}G_{-\delta}(X_s,Y_s)\dd s\Big\}\Big]\cr
 \ar=\ar
\mbf{E}\Big[
(X_T^{-\beta}Y_T)^{-\delta}\exp\Big\{\delta\int_0^TG_{-\delta}(X_s,Y_s)\dd s\Big\}\Big].
 \eeqnn
By \eqref{3.2} and \eqref{3.4}, for $s<\tau_w$ and $\sup_{s\ge0}(X_s\vee Y_s)\le \varepsilon_2$, we have $G_{-\delta}(X_s,Y_s)\ge0$.
Then
 \beqnn
1
 \ar\ge\ar
\mbf{E}\Big[
(X_T^{-\beta}Y_T)^{-\delta}\exp\Big\{\delta\int_0^TG_{-\delta}(X_s,Y_s)\dd s\Big\}1_{\{\sup_{s\ge0}(X_s\vee Y_s)\le
\varepsilon_2,\tau_w<\sigma_v\}}\Big] \cr
 \ar\ge\ar
w^{-\delta}
\mbf{P}\Big\{\sup_{s\ge0}(X_s\vee Y_s)\le \varepsilon_2,\,\tau_w<\sigma_v\Big\}.
 \eeqnn
By Lemma \ref{t2.1} (i), we have
$\lim_{v\to\infty}\sigma_v=\infty$ almost surely.
Letting $v\to\infty$ in the above inequality we obtain
 \beqnn
\mbf{P}\Big\{\sup_{s\ge0}(X_s\vee Y_s)\le \varepsilon_2,\,\tau_w<\infty\Big\}
\le w^\delta.
 \eeqnn
Combining  \eqref{3.3} and \eqref{3.4} it follows that
 \beqnn
 \ar\ar
\mbf{P}\{\tau_w=\infty\}
=1-\mbf{P}\{\tau_w<\infty\} \cr
 \ar\ge\ar
1-\mbf{P}\Big\{\sup_{s\ge0}(X_s\vee Y_s)\le \varepsilon_2,\,\tau_w<\infty\Big\}
-\mbf{P}\Big\{\sup_{s\ge0}(X_s\vee Y_s)\ge \varepsilon_2\Big\} \cr
 \ar\ge\ar
1-w^\delta
-C[\varepsilon_2^{\beta/4}
+\varepsilon_2^{1/(4\beta)}]>0,
 \eeqnn
which implies
$\mbf{P}\{\tau_0^Y=\infty\}>0$
and ends the proof.
\qed

\noindent{\it Proof of Theorem \ref{t0.05}}.
By Lemma \ref{t4.6},
without loss of generality we assume that $X_0<\varepsilon$ and $\mbf{P}\big\{ \sup_{t\geq 0}X_t\leq \varepsilon\big\}>0$.
Let $(\bar{Y}_t)_{t\ge0}$ be the solution to \eqref{3.18}
with $C_n$ replaced by $C_\varepsilon:=\sup_{x\in[0,\varepsilon]}\kappa(x)$.
By the comparison theorem (Proposition \ref{t4.4}), we have
 \beqnn
\mbf{P}\Big\{ Y_t\geq \bar{Y}_t \,\,\,\text{for all}\,\,\, t\geq 0\big|\sup_{t\geq 0}X_t\leq \varepsilon \Big\}=1.
 \eeqnn
Since $\mbf{P}\{\sigma^{\bar{Y}}(\varepsilon^\beta)<\infty \}>0$
by \cite[Proposition 2.11]{LYZh}
and $\bar{Y}$ is independent of $X$, then
 \beqnn
 \ar\ar
\mbf{P}\big\{X(\sigma^{\bar{Y}}(\varepsilon^\beta)) \leq \varepsilon, Y(\sigma^{\bar{Y}}(\varepsilon^\beta))\ge \varepsilon^\beta,
\sigma^{\bar{Y}}(\varepsilon^\beta)<\infty\big\} \cr
 \ar\ge\ar
\mbf{P}\big\{X(\sigma^{\bar{Y}}(\varepsilon^\beta)) \leq \varepsilon, \bar{Y}(\sigma^{\bar{Y}}(\varepsilon^\beta))\ge \varepsilon^\beta,
\sigma^{\bar{Y}}(\varepsilon^\beta)<\infty,\sup_{t\geq 0}X_t\leq \varepsilon\big\} \cr
 \ar=\ar
\mbf{P}\big\{\sup_{t\geq 0}X_t\leq \varepsilon,\sigma^{\bar{Y}}(\varepsilon^\beta)<\infty\big\}
=
\mbf{P}\big\{ \sup_{t\geq 0}X_t\leq \varepsilon\big\} \mbf{P}\big\{\sigma^{\bar{Y}}(\varepsilon^\beta)<\infty\big\}
>0.
 \eeqnn

Note that by Proposition \ref{t4.4} and Lemma \ref{t3.3},
there exist constants $\beta>0$
and small $\varepsilon>0$ so that
$\mbf{P}\{\tau_0^Y=\infty\}>0$
if $X_0\le\varepsilon,Y_0\ge \varepsilon^\beta$. Applying the strong Markov property to process $(X, Y)$ at time $\sigma^{\bar{Y}}(\varepsilon^\beta)$,
we have
$\mbf{P}\{\tau_0^Y=\infty\}>0$ for any $Y_0>0$, and the proof is completed.
\qed

\bremark
Under Condition \ref{c} (ii) and (iii),
for $p>0$, the corresponding estimate of $G_0(x,y)$ in Step 1 of the proof of Lemma \ref{t3.3} is not easy to  establish and thus
the approach of showing Theorem \ref{t0.05} does not appear to be valid for the proofs of Theorem \ref{t0.02}
and Conjecture \ref{conje}.
\eremark

\section{Appendix}
\setcounter{equation}{0}

In this section we present a proof of the comparison theorem in Proposition  \ref{t4.4}
and show in Lemma
\ref{t5.1} that if all the functions in \eqref{1.1} are locally Lipschitz, then
\eqref{1.1} has a pathwise unique solution.
Consequently, $(X,Y)$ is a Markov process.

\noindent{\it Proof of Proposition \ref{t4.4}.}
For $k\ge1$ define
 \beqnn
h_k:=\exp\{-k(k+1)/2\}.
 \eeqnn
Let $\psi_k$ be a nonnegative function on $\mbb{R}$
with support in $(h_k,h_{k-1})$,
$\int_{h_k}^{h_{k-1}}\psi_k(x)\dd x=1$
and
 \beqnn
0\le\psi_k(x)\le
2k^{-1}x^{-1}1_{(h_k,h_{k-1})}(x).
 \eeqnn
For $x\in\mbb{R}$ and $k\ge1$ let
 \beqnn
\phi_k(x):=1_{\{x>0\}}\int_0^{x}\dd y\int_0^y\psi_k(z)\dd z.
 \eeqnn
For $k\ge1$ and $y,z\in\mbb{R}$ put
 \beqlb\label{DN}
\mathcal{D}_k(y,z):=\phi_k(y+z)-\phi_k(y)-z\phi_k'(y).
 \eeqlb
For $t\ge0$ let $\bar{x}(t)=x_1(t)-x_2(t)$,
$\hat{B}(t)=B_1(t,x_2(t))-B_2(t,x_2(t))$,
$\bar{B}(t)=B_1(t,x_1(t))-B_1(t,x_2(t))$,
$\bar{U}(t)=U(x_1(t))-U(x_2(t))$ and
$\bar{V}(t)=V(x_1(t))-V(x_2(t))$.
It then follows from \eqref{8.1} that
 \beqnn
\bar{x}(t\wedge\tilde{\gamma}_n)
 \ar=\ar
\bar{x}(0)
+\int_0^{t\wedge\tilde{\gamma}_n}[\hat{B}(s)+\bar{B}(s)]\dd s
+\int_0^{t\wedge\tilde{\gamma}_n} \bar{U}(s)\dd W_s \cr
 \ar\ar
+\int_0^{t\wedge\tilde{\gamma}_n}\int_0^\infty\int_0^\infty
g(s-,u)z \tilde{N}(\dd s,\dd z,\dd u),
 \eeqnn
where $\tilde{\gamma}_n$ is defined in \eqref{4.2} and $g(s,u):=1_{\{u\le V(x_1(s))\}}-
1_{\{u\le V(x_2(s))\}}$.
Using It\^o's formula we obtain
 \beqnn
 \ar\ar
\phi_k(\bar{x}(t\wedge\tilde{\gamma}_n)) \cr
 \ar=\ar
\phi_k(\bar{x}(0))
+\int_0^{t\wedge\tilde{\gamma}_n}\phi_k'(\bar{x}(s))[\hat{B}(s)+\bar{B}(s)]\dd s
+\frac12\int_0^{t\wedge\tilde{\gamma}_n}\phi_k''(\bar{x}(s))\bar{U}(s)^2\dd s \cr
 \ar\ar
+\int_0^{t\wedge\tilde{\gamma}_n}\dd s\int_0^\infty \bar{V}(s)\mathcal{D}_k(\bar{x}(s),\mbox{sgn}(\bar{V}(s))z )\nu(\dd z)
+\mbox{mart.},
 \eeqnn
where $\mbox{sgn}(x)=1_{\{x>0\}}-1_{\{x<0\}}$.
It follows that
 \beqlb\label{8.2}
 \ar\ar
\mbf{E}\big[\phi_k(\bar{x}(t\wedge\tilde{\gamma}_n))\big] \cr
 \ar=\ar
\phi_k(\bar{x}(0))
+
\mbf{E}\Big[\int_0^{t\wedge\tilde{\gamma}_n}\phi_k'(\bar{x}(s))\hat{B}(s)\dd s\Big] \cr
 \ar\ar
+\mbf{E}\Big[\int_0^{t\wedge\tilde{\gamma}_n}\phi_k'(\bar{x}(s))\bar{B}(s)\dd s\Big]
+
\frac12\,\mbf{E}\Big[\int_0^{t\wedge\tilde{\gamma}_n}
\phi_k''(\bar{x}(s))\bar{U}(s)^2\dd s\Big] \cr
 \ar\ar
+\mbf{E}\Big[\int_0^{t\wedge\tilde{\gamma}_n}\dd s\int_0^\infty \bar{V}(s)\mathcal{D}_k(\bar{x}(s),\mbox{sgn}(\bar{V}(s))z )\nu(\dd z)\Big] \cr
 \ar=:\ar
\phi_k(\bar{x}(0))
+\sum_{i=1}^4I_{n,k}^i(t).
 \eeqlb
For $x\in\mbb{R}$, $x^+:=x\vee0$.
By \cite[Lemma 2.1]{XiongYang17},
\[\lim_{k\to\infty}\phi_k(x)=x^+,\quad
\lim_{k\to\infty}\phi_k'(x)=1_{\{x>0\}},\quad
|x|\phi_k''(x)\le 2k^{-1},\]
and
\[\mathcal{D}_k(y,z)\le (2k^{-1}z^2/y)\wedge(2|z|)\quad
\text{for all}~\, k\ge1, x,y\in\mbb{R} \,\,\text{and}\,\, z\ge0\,\,\text{ with}\,\, y(y+z)>0.\]
Then using the assumptions and the dominate convergence,
 \beqnn
\lim_{k\to\infty}\mbf{E}\big[\phi_k(\bar{x}(t\wedge\tilde{\gamma}_n))\big]
=\mbf{E}\big[\bar{x}^+(t\wedge\tilde{\gamma}_n)\big],\quad
\lim_{k\to\infty}I_{n,k}^1(t)
=\mbf{E}\Big[\int_0^{t\wedge\tilde{\gamma}_n}1_{\{\bar{x}(s)>0\}}\hat{B}(s)\dd s\Big]
\le0
 \eeqnn
and
 \beqnn
\lim_{k\to\infty}I_{n,k}^2(t)
\le
C_n\mbf{E}\Big[\int_0^{t\wedge\tilde{\gamma}_n}\bar{x}^+(s)\dd s\Big],~
\lim_{k\to\infty}I_{n,k}^3(t)
=
\lim_{k\to\infty}I_{n,k}^4(t)=0.
 \eeqnn
Combining with \eqref{8.2} we get
 \beqnn
\mbf{E}\big[\bar{x}^+(t\wedge\tilde{\gamma}_n)\big]
\le
\bar{x}^+(0)+C_n\mbf{E}\Big[\int_0^{t\wedge\tilde{\gamma}_n}\bar{x}^+(s)\dd s\Big]
\le
C_n\int_0^t\mbf{E}\big[\bar{x}^+(s\wedge\tilde{\gamma}_n)\big]\dd s.
 \eeqnn
From Gronwall's lemma it follows that $\mbf{E}\big[\bar{x}^+(t\wedge\tilde{\gamma}_n)\big]=0$.
Letting $n\to\infty$ we get
$\bar{x}^+(t\wedge\tilde{\gamma})=0$ almost surely for each
fixed $t>0$.
By the right continuity of $t\mapsto x_i(t)$ $(i=1,2)$ we concludes the
proof.
\qed

\blemma\label{t5.1}
Suppose that the functions
$a_i,b_i, i=1,2,3$ and $\theta,\kappa$ in \eqref{1.1} are all locally Lipschitz, i.e.,
for each $m,n\ge1$,  there is a constant $C_{m,n}>0$ so that
 \beqnn
\sum_{i=1,2,3}\big[|a_i(x)-a_i(y)|+|b_i(x)-b_i(y)|\big]\le
C_{m,n}|x-y|,
\qquad x,y\in[n^{-1}, m].
 \eeqnn
Then
SDE \eqref{1.1} has a nonnegative pathwise unique solution.
\elemma
\proof
For $n\ge1$ and $i=1,2,3$ let $a_i^n(x):=a_i((x\wedge n)\vee n^{-1})$.
Define $b_i^n$ and $\theta^n,\kappa^n,$
similarly.
Inspire by the argument in \cite[Theorem 3.1]{HLX},
let
 \beqnn
U=\{1,2\}\times(0,\infty)^2,\quad
U_0=\big(\{1\}\times(0,1)\times(0,\infty)\big)
\cup\big(\{2\}\times(0,1)\times(0,\infty)\big).
 \eeqnn
Let
 \beqnn
\tilde{N}_0(\dd s,\dd v,\dd z,\dd u)
:=
\delta_1(\dd v)\tilde{M}(\dd s,\dd z,\dd u)
+\delta_2(\dd v)\tilde{N}(\dd s,\dd z,\dd u).
 \eeqnn
Then $\tilde{N}_0$ is a compensated Poisson random measure on
$(0,\infty)\times U$ with intensity
 \beqnn
\dd s[\delta_1(\dd v)\mu(\dd z)+\delta_2(\dd v)\nu(\dd z)]\dd u.
 \eeqnn
Define functions $f_1^n$ and $f_2^n$ on $[0,\infty)\times U$ by
 \beqnn
f_1^n(x,v,z,u):=
z 1_{\{v=1,u\le a_3^n(x)\}},
\quad
f_2^n(x,v,z,u):=
z 1_{\{v=2,u\le b_3^n(x)\}}.
 \eeqnn
Write $\textbf{u}$ for $(v,z,u)$.
Then SDE \eqref{1.1} can be written into this form
 \beqlb\label{1.11}
  \left\{
   \begin{aligned}
   X_t &=X_0 -\int_0^t \Big[a_1^n(X_s)+a_3^n(X_s)\int_1^\infty z\mu(\dd z)\Big] \dd s+\int_0^t a_2^n(X_s)^{1/2}\dd B_s\\ &\quad
   +\int_0^t\int_{U_0} f_1^n(X_{s-},\textbf{u})\tilde{N}_0(\dd s,\dd \textbf{u})
   +\int_0^t\int_{U\setminus U_0} f_1^n(X_{s-},\textbf{u})N_0(\dd s,\dd \textbf{u}), \\
   Y_t &=Y_0+\int_0^t\Big[-b_1^n(Y_s)-\theta^n(Y_s)\kappa^n(X_s)
+b_3^n(X_s)\int_1^\infty z\nu(\dd z)\Big]\dd s
   +\int_0^t b_2^n(Y_s)^{1/2}\dd W_s\\ &\quad
   +\int_0^t\int_{U_0} f_2^n(Y_{s-},\textbf{u})\tilde{N}_0(\dd s,\dd \textbf{u})
+\int_0^t\int_{U\setminus U_0} f_2^n(Y_{s-},\textbf{u})N_0(\dd s,\dd \textbf{u}),
   \end{aligned}
   \right.
 \eeqlb
where $N_0$ is the corresponding Poisson random measure of $\tilde{N}_0$.
It follows from \cite[p. 245]{IkWa} that \eqref{1.11} has a strong unique
solution $(X_t^n,Y_t^n)_{t\ge0}$.
Letting $n\to\infty$,
by the same argument in \cite[Theorem 3.1 (i)]{LYZh},
SDE \eqref{1.1} has a pathwise unique solution.

\section*{Acknowledgments}
We thank anonymous referees and Rugang Ma for very detailed and helpful comments.
The last two authors thank the support from a
Foreign Experts Project of North Minzu University.

\end{document}